\documentclass[12pt]{amsart} 

\usepackage[dvips]{graphics}

\usepackage{epsfig}
\usepackage{rotating}

\usepackage{amsfonts}
\usepackage{amsmath, amsthm, amssymb, ulem, amscd} 

\def\crn#1#2{{\vcenter{\vbox{
        \hbox{\kern#2pt \vrule width.#2pt height#1pt
           }
          \hrule height.#2pt}}}}
\def\intprod{\mathchoice\crn54\crn54\crn{3.75}3\crn{2.5}2}
\def\into{\mathbin{\intprod}}

\newcommand{\congvert}
{\begin{turn}{-90}
$\cong\text{ }$ 
\end{turn}}

\newcommand{\eqvert}
{\begin{turn}{-90}
$=\text{ }$ 
\end{turn}}

\newcommand{\pa}{\partial}

\newcommand{\st}{{\tilde{s}}}
\newcommand{\lt}{{\tilde{l}}}

\newcommand{\wt}{{\tilde{\bigwedge}}}
\newcommand{\Lt}{{\tilde{\textstyle{\bigwedge}}}} 

\newcommand{\trt}{{\tilde{\operatorname{tr}}}}

\newcommand{\cG}{{\mathcal G}}
\newcommand{\cA}{{\mathcal A}}

\newcommand{\cE}{{\mathcal E}}
\newcommand{\cM}{{\mathcal M}}
\newcommand{\cGt}{{\tilde {\mathcal G}}}
\newcommand{\cRt}{{\tilde {\mathcal R}}}

\newcommand{\cTt}{{\tilde {\mathcal T}}}

\newcommand{\cEt}{{\tilde {\mathcal E}}}
\newcommand{\cJt}{{\tilde {\mathcal J}}}

\newcommand{\cO}{{\mathcal O}}

\newcommand{\cR}{{\mathcal R}}
\newcommand{\cJ}{{\mathcal J}}
\newcommand{\cD}{{\mathcal D}} 
\newcommand{\cQ}{{\mathcal Q}}
\newcommand{\cV}{{\mathcal V}}
\newcommand{\cU}{{\mathcal U}}

\newcommand{\cN}{{\mathcal N}}

\newcommand{\gh}{{\widehat g}}
\newcommand{\Ah}{{\widehat A}} 
\newcommand{\gt}{{\tilde g}}
\newcommand{\htt}{{\tilde h}}

\newcommand{\ft}{{\tilde f}}
\newcommand{\Rt}{{\tilde R}}
\newcommand{\cH}{{\mathcal H}}
\newcommand{\Z}{{\mathbb Z}}
\newcommand{\R}{{\mathbb R}}

\newcommand{\N}{{\mathbb N}}

\newcommand{\C}{{\mathbb C}}

\newcommand{\dt}{\tilde \delta} 
\newcommand{\Dt}{\tilde \Delta} 

\newcommand{\Ric}{\operatorname{Ric}}

\newcommand{\tf}{\operatorname{tf}}

\renewcommand{\tilde}{\widetilde}

\newcommand{\Diff}{\text{\rm Diff}}
\newcommand{\CDiff}{\text{\rm CDiff}}
\newcommand{\ODiff}{\text{\rm ODiff}}

\theoremstyle{plain}
\newtheorem{theorem}{Theorem}[section]
\newtheorem{lemma}[theorem]{Lemma}
\newtheorem{proposition}[theorem]{Proposition}

\theoremstyle{definition}

\newtheorem{definition}[theorem]{Definition}

\theoremstyle{remark}

\numberwithin{equation}{section}

\title{Jet Isomorphism for Conformal Geometry}  

\author{C. Robin Graham}
\address{Department of Mathematics, University of Washington,
Box 354350\\
Seattle, WA 98195-4350}
\email{robin@math.washington.edu}

\begin{document}

\maketitle

\thispagestyle{empty}

\renewcommand{\thefootnote}{}
\footnotetext{This research was partially supported by NSF grant \# DMS 
 0505701.}

\section*{Introduction}\label{intro}

Local invariants of a metric in Riemannian geometry are quantities
expressible in local coordinates in terms of the metric and its derivatives  
and which have an invariance property under changes of coordinates.        
It is a fundamental fact that such invariants may be written in terms of
the 
curvature tensor of the metric and its covariant derivatives.  In this
form, they can be identified with invariants of the orthogonal group acting   
algebraically on the space of 
possible curvature tensors and derivatives.  We refer to the result   
asserting that the space of infinite order jets of metrics modulo
coordinate changes is isomorphic to a space of curvature 
tensors and derivatives modulo the orthogonal group as a jet isomorphism
theorem.  Such results recast the study and
classification of local geometric invariants in purely algebraic terms,
in which form the methods of invariant theory and representation
theory can be brought to bear.  

The goal of this paper is to describe analogous jet isomorphism theorems
in the context of conformal geometry.  In conformal geometry one is given a
metric only up to scale.  The results in the conformal case provide a
tensorial description of the space of jets of metrics   
modulo changes of coordinates and conformal factor.  The motion group of
the flat model is the conformal group $G=O(n+1,1)/\{\pm I\}$ acting
projectively on the sphere $S^n$ and the role   
of the orthogonal group in Riemannian geometry is played by the parabolic
subgroup $P\subset G$ 
preserving a null line.  Since $P$ is a matrix group in $n+2$ dimensions,
its natural tensor representations are on tensor powers of $\R^{n+2}$.
Thus one expects the appearance of tensors in $n+2$ dimensions in
conformal jet isomorphism theorems.  

When $n$ is odd, the ambient metric construction of \cite{FG1} associates
to a conformal 
Riemannian manifold $(M,[g])$ of dimension $n$ an infinite order jet 
of a Lorentzian metric $\gt$ along a hypersurface 
in a space $\cGt$ of dimension $n+2$, uniquely determined up to
diffeomorphism.  The tensors in the odd-dimensional conformal jet
isomorphism   
theorem are the curvature tensor and its covariant derivatives for the
ambient metric.  They satisfy extra identities beyond those satisfied 
by the derivatives of curvature of a general metric as a consequence of the
Ricci-flatness and homogeneity conditions satisfied by an ambient metric.   
The elaboration of these identities and a formulation and proof of a 
conformal jet isomorphism theorem in odd dimensions are given in
\cite{FG2}.  The algebra of the proof is more involved than in the
Riemannian case because one is comparing tensors in different dimensions.     

In even dimensions, the ambient
metric construction is obstructed at finite order, so this gives only a  
finite order version of a jet isomorphism theorem.  
An extension of the ambient metric construction to all orders in even
dimensions has recently been described in \cite{GH1}.  Based on this, joint 
work in preparation with K. Hirachi formulates and   
proves an infinite order version of the jet isomorphism theorem in even
dimensions.  The 
method of proof of the jet isomorphism theorem in this work is different
than that used in \cite{FG2}; it relies on an ambient lift of the conformal  
deformation complex on $G/P$.  This same method also can be used   
to give another proof of the odd-dimensional jet isomorphism theorem.  
This proof in the odd-dimensional case will be sketched in \S\ref{deform}
below and the details will be the subject of \cite{GH2}.     

In the
even-dimensional case, an ambient metric depends not only on the conformal 
manifold $(M,[g])$, but also on the choice of a trace-free symmetric
2-tensor called the ambiguity tensor.  Likewise, the even-dimensional jet
isomorphism theorem provides a tensorial description of an enlargement of
the space of jets of metrics by the space of jets of the ambiguity.

The results concerning the ambient lift of the deformation complex are
perhaps of independent interest.  In odd dimensions, all the spaces
occuring in the complex except for the next to last have isomorphic   
realizations in terms of infinite order jets along a hypersurface of  
tensors defined on the ambient space, and the maps in the complex simplify
when written in these realizations.  The situation in even 
dimensions is more complicated owing to the existence of
ambiguities in the lifts, but it is nonetheless possible to prove results
concerning ambient realization including ambiguities which can be used to    
prove the jet isomorphism theorem in even dimensions.  

In \S\ref{odd}, we first recall 
the jet isomorphism theorem for pseudo-Riemannian geometry.
We then show how the space of jets of metrics   
modulo changes of coordinates and conformal factor has a natural action of 
$P$.   In both cases, our presentation is in terms of a quotient of the   
space of jets of all metrics as in \cite{GH2} rather than   
in terms of metrics in geodesic normal coordinates as in 
\cite{FG2}.  Next we   
formulate the odd-dimensional conformal jet isomorphism theorem.
We also briefly review the ambient metric construction 
and show how in odd dimensions 
it gives rise to the map from jets of metrics to the algebraic   
space of ambient curvature tensors and their covariant derivatives.   
In \S\ref{deform}, we begin by introducing the conformal deformation
complex.  We then formulate Theorem~\ref{amlift}, which gives the ambient
realization 
in the unobstructed cases for sections of homogeneous bundles on $G/P$ with
symmetries defined by a Young diagram with no more than two columns.  
This generalizes to tensors results of \cite{EG} for scalars.  We sketch 
the proof for scalars and 1-forms.  A consequence of Theorem~\ref{amlift}
is the ambient realization of the deformation complex 
for $n$ odd.  We conclude \S\ref{deform} by sketching the proof of 
the odd-dimensional jet isomorphism theorem, 
using this ambient realization and the exactness of the deformation complex
on jets.   
In \S\ref{even}, we first discuss the ambient lift in the simplest
obstructed case: that of densities whose weight is such that smooth
harmonic extension is obstructed.  We indicate how an infinite order
harmonic extension 
involving a log term can always be found, albeit with an ambiguity, and how
to reformulate the harmonic extension in terms of the smooth part.  These
considerations result in Theorem~\ref{ambscalar}, the substitute ambient  
lift theorem for scalars in the obstructed cases.   The discussion of the
scalar case illustrates the phenomena which occur for jets of conformal
structures in even dimensions.  We next formulate the   
conformal jet isomorphism theorem in even dimensions.  Then we outline the 
extension of the ambient metric construction to all orders and indicate 
by analogy with \S\ref{odd} how it gives rise to the map from jets of
metrics and ambiguity tensors to ambient curvature tensors and briefly
indicate what is involved in the   
proof of the even-dimensional jet isomorphism theorem.  

A jet isomorphism theorem for a parabolic geometry was first considered in  
\cite{F}, and the entire perspective explicated here owes much to this
pioneering work.  The idea of incorporating an ambiguity into
a jet isomorphism theorem was introduced in \cite{H}.  A different approach
for conformal geometry using tractor calculus rather than the ambient metric to 
construct curvature tensors is given in \cite{G}.

The lectures on which this paper is based were delivered at the 2007  Winter
School 'Geometry and Physics' at Srn\'i.  The author is grateful to the
organizing committee, particularly  to Vladimir Sou\u{c}ek, for the
opportunity to participate in this school.

\section{Jet isomorphism, odd dimensions}\label{odd}

We begin by reviewing the jet isomorphism theorem for pseudo-Riemannian
geometry.  Fix a signature $(p,q)$, $p+q=n\geq 2$, and a reference
quadratic form $h_{ij}$ on $\R^n$ of signature $(p,q)$  
(one
typically takes $h_{ij}=\delta_{ij}$ in the positive definite case).  
By a change of coordinates, any metric of signature $(p,q)$ can be made to
equal $h_{ij}$ at the origin.  It is convenient to include this 
normalization in our definition.  So we set
$$
\cM =\{\text{Jets of metrics $g_{ij}$ such that $g_{ij}(0)=h_{ij}$}\}.
$$
Here jets means infinite order jets of smooth
metrics at the origin in $\R^n$. 
We can identify an element of $\cM$ with the list 
$(\pa^{\alpha}g_{ij}(0))_{|\alpha|\geq 1}$.  Define also 
$$
\Diff =\{\text{Jets of local diffeomorphisms $\varphi$ of $\R^n$
satisfying $\varphi(0)=0$}\}.
$$
Then $\Diff$ is a group under composition.  
Since we have normalized our metric at the origin, we need to restrict to  
diffeomorphisms which preserve the normalization.  So we define the
subgroup $\ODiff\subset \Diff $ by
$$
\ODiff =
 \{\varphi \in \Diff: \varphi'(0)\in O(h)\} 
$$
and the normal subgroup $\Diff_0\subset \ODiff $ by 
$$
\Diff_0 =
 \{\varphi: \varphi'(0)=I\}.
$$
Then $\ODiff$ acts on $\cM$ on the left by
$\varphi.g=(\varphi^{-1})^*g$.  We can view $O(h)$ as the subgroup of 
$\ODiff$ consisting of linear transformations.  Then $O(h)$ is the isotropy
group in $\ODiff$ of the flat metric $h\in \cM$, and $\ODiff = O(h)\cdot  
\Diff_0$.  

Since $\Diff_0$ is a normal subgroup of $\ODiff$, there is an induced
action of $O(h)$ on the orbit space $\cM/\Diff_0$ (we write the quotient on
the right even though this is a left action).  Local invariants of
pseudo-Riemannian 
metrics can be thought of as functions on $\cM$ which are invariant under
the action of $\Diff_0$ and equivariant under $O(h)$;  such a function
determines an assignment to each metric on an arbitrary manifold of a
section of the associated bundle by evaluation at each 
point in local coordinates.  The jet
isomorphism theorem for pseudo-Riemannian geometry provides an
$O(h)$-equivariant description of the space $\cM/\Diff_0$ in terms of
curvature tensors and their covariant derivatives.  

\begin{definition}\label{R}
The space $\cR\subset \prod_{r=0}^\infty
\bigwedge^2\R^n{}^*\otimes\bigwedge^2\R^n{}^*\otimes \bigotimes^r\R^n{}^*$
is the set of lists $(R^{(0)},R^{(1)},R^{(2)},\cdots )$ with
$R^{(r)}\in\bigwedge^2\R^n{}^*\otimes\bigwedge^2\R^n{}^*\otimes 
\bigotimes^r\R^n{}^*$, such that:  
\begin{enumerate}
\item
$R_{i[jkl],m_1\cdots m_r}=0$
\item
$R_{ij[kl,m_1]m_2\cdots m_r} = 0$
\item
$R_{ijkl,m_1\cdots [m_{s-1}m_s]\cdots m_r}= Q^{(s)}_{ijklm_1\cdots
   m_r}(R)$.
\end{enumerate}
Here the comma after the first four indices is just a marker separating
these indices.  
$Q^{(s)}_{ijklm_1\cdots m_r}(R)$ denotes the quadratic
expression in the $R^{(r')}$ with $r'\leq r-2$ which one
obtains by covariantly differentiating the usual Ricci identity for
commuting covariant  
derivatives, expanding the differentiations using the Leibnitz rule, and
then setting equal
to $h$ the metric which contracts the two factors in each term.     
We have suppressed the ${}^{(r)}$ on the $R^{(r)}$ since the value of $r$
is evident from the list of indices.  
\end{definition}  

The action of $O(h)$ on $\R^n$ induces actions on the spaces of tensors in
the usual way and therefore also on 
$\prod_{r=0}^\infty \bigwedge^2\R^n{}^*\otimes\bigwedge^2\R^n{}^*\otimes
\bigotimes^r\R^n{}^*$.  Since $\cR$ is an $O(h)$-invariant subset of this
product, $\cR$ has a natural $O(h)$ action.   

Evaluation of the covariant derivatives of curvature of a metric at
the origin induces a polynomial map $\cM\rightarrow \cR$.  Since the
covariant derivatives of curvature are tensors, it follows that this map
passes to the quotient, and so defines a map $\cM/\Diff_0\rightarrow \cR$
which is $O(h)$-equivariant.  The pseudo-Riemannian jet isomorphism theorem
is then the following.

\begin{theorem}
The map $\cM/\Diff_0\rightarrow \cR$ is an $O(h)$-equivariant bijection
with polynomial inverse.
\end{theorem}

\noindent
The proof proceeds via the introduction of geodesic normal coordinates.  
These provide a slice for the action of $\Diff_0$ on 
$\cM$:  each orbit in $\cM/\Diff_0$ is represented by a unique jet of a
metric for which the background coordinates on $\R^n$ are geodesic normal
coordinates to infinite order at the origin.  A linearization argument
reduces the theorem to showing that the linearized map restricted to  
infinitesimal jets of metrics in normal coordinates is a vector space
isomorphism.  The linearized   
map can be explicitly identified as the direct sum over $r$ of intertwining
maps between two equivalent realizations corresponding to different Young
projectors of irreducible representations of $GL(n,\R)$.  See \cite{E} for
the analysis of the linearized map in a similar context.   

This jet isomorphism theorem is fundamental in consideration of local
pseudo-Riemannian invariants.  It shows that such invariants correspond
exactly to 
$O(h)$-invariants of $\cR$.  Weyl's classical invariant theory for $O(h)$
completely describes such invariants.  

Next we pass to the conformal analogue.  We begin with a discussion of the
conformal 
group and its parabolic subgroup $P$ which plays the role in conformal
geometry of the group $O(h)$ in pseudo-Riemannian geometry.  In the
conformal case we assume that $n\geq 3$.  

Define a quadratic form $\htt$ of signature $(p+1,q+1)$ on $\R^{n+2}$ by 
$$
\htt_{IJ}=
\begin{pmatrix}
0&0&1\\
0&h_{ij}&0\\
1&0&0
\end{pmatrix}.
$$
On $\R^{n+2}$ we use as coordinates $x^I=(x^0,x^i,x^\infty)$.   
The null cone of $\htt$ is
$$
\cN = \{x\in \R^{n+2}\setminus \{0\}: 
\htt_{IJ}x^Ix^J =0\},
$$
whose projectivization is the quadric
$$
\cQ = \{[x]\in \mathbb{P}^{n+1}:x\in \cN\},
$$
with projection $\pi:\cN\rightarrow \cQ$.  
If $x\in \cN$, the metric 
$\htt_{IJ}dx^Idx^J$ on $\R^{n+2}$ is degenerate when restricted 
to $T_x\cN$:  $\htt(X,Y)=0$ for all $Y\in T_x\cN$,
where $X=x^I\pa_I$ is the Euler field.  Consequently,   
$\htt|_{T_x\mathcal{N}}$ determines a nondegenerate quadratic form on
$T_x\mathcal{N}/\operatorname{span}{X}\cong T_{\pi(x)}\cQ$.  
As $x$ varies over a line in $\mathcal{N}$, the  
resulting quadratic forms on $T_{\pi(x)}\cQ$ define a metric up to scale,
i.e. a conformal class of metrics on $\cQ$ of signature $(p,q)$.

The conformal group is $G=O(\htt)/\{\pm I\}$.   The linear 
action of 
$O(\htt)$ on $\R^{n+2}$ preserves $\cN$, so there is an induced
action of  
$G$ on $\cQ$ which is transitive.  Since $O(\htt)$ acts by isometries of
the metric $\htt_{IJ}dx^Idx^J$ on $\R^{n+2}$, the induced action of $G$ on
$\cQ$ is by conformal transformations.  Define the subgroup $P\subset G$ to
be the isotropy group of $[e_0]\in \cQ$, so that $\cQ=G/P$.  Then $P$ can
be identified with the  
subgroup $P=\{p\in O(\htt): pe_0 = a e_0,\, a>0\}$.  The first column of $p$
is $ae_0$; combining this with the fact that $p\in O(\htt)$, one finds
that $p\in P$ is of the form 
\begin{equation}\label{pform}
p=\left(
\begin{matrix}
a&b_j&c\\
0&m^i{}_j&d^i\\
0&0&a^{-1}
\end{matrix}
\right),
\end{equation}
where 
$$
a>0,\,\, b_j\in \R^n{}^*,\,\, m^i{}_j\in O(h)\quad \text{and}\quad
c=-\frac{1}{2a}b_jb^j,\,\, d^i =-\frac{1}{a}m^{ij}b_j.
$$ 
Lower case indices are raised and lowered using $h$.  
The parameters $a>0$, $b_j\in \R^n{}^*$ and $m^i{}_j\in O(h)$ are free, so
that $P$ can be written as the product of its subgroups:
$$
P=\mathbb{R}_+\cdot \R^n \cdot O(h).
$$

Since $\htt_{IJ}x^Ix^J=2x^0x^\infty +|x|^2$, where $|x|^2= h_{ij}x^ix^j$,
the intersection of $\cQ$ with 
the cell $\{[x^I]: x^0\neq 0\}$ can be identified with $\R^n$ via 
the inclusion $i:\R^n\to \cQ$ defined by 
$$
i(x)=
\left[
\begin{matrix}
1\\
x\\  
-\frac12 |x|^2
\end{matrix}
\right].
$$
In this identification, the conformal structure on $\cQ$ is represented by 
the flat metric $h_{ij}dx^idx^j$ on $\R^{n}$.  For $p\in P$, we will denote
by $\varphi_p$ the corresponding conformal transformation on $\R^n$, and by
$\Omega_p$ the conformal factor, so that 
$$
\varphi_p^*h =\Omega_p^2\, h.
$$   
These are given explicitly by:
$$
(\varphi_p(x))^i= \frac{m^i{}_jx^j -\frac12 |x|^2 d^i}{a+b_jx^j-\frac12 c
|x|^2},\qquad\quad
\Omega_p = (a+b_jx^j-\tfrac12 c|x|^2)^{-1}. 
$$
Observe that 
\begin{equation}\label{initial}
\varphi_p'(0)=a^{-1}m^i{}_j,\qquad
\Omega_p(0)= a^{-1},\qquad \Omega_p'(0)=-a^{-2}b_j.   
\end{equation}

We now consider changing the metric by rescaling 
as well as by diffeomorphism.  Set 
$$
C^{\infty}_+ =\{\text{Jets at $0\in \R^n$ of smooth positive 
functions $\Omega$}\}.  
$$
Then $C^\infty_+$ is a group under multiplication.  
Consider the semidirect product group $\Diff\ltimes C^\infty_+$, where the
product is defined so that $g\rightarrow (\varphi^{-1})^*(\Omega^2g)$ 
is an action.  This product is given explicitly by:  
$$
(\varphi_1,\Omega_1)\cdot
(\varphi_2,\Omega_2)=(\varphi_1\circ\varphi_2,(\Omega_1\circ 
\varphi_2) \,\Omega_2).
$$
As before, we need to preserve the normalization $g_{ij}(0)=h_{ij}$.  So we 
define the subgroup
$\CDiff\subset \Diff\ltimes C^\infty_+$ 
by
$$
\CDiff  = \{(\varphi,\Omega): (\Omega^{-1}\varphi')(0)\in O(h)\}
$$
and we set 
$$
\CDiff_0 = \{(\varphi,\Omega): \varphi'(0)
=I,\quad\Omega(0) =1,\quad d\Omega(0)=0\}\subset \CDiff.    
$$
Then
$\CDiff$ acts on $\cM$ by:  $(\varphi,\Omega).g=
(\varphi^{-1})^*(\Omega^2g)$.  
We can view $P\subset \CDiff$ by $p\mapsto (\varphi_p,\Omega_p)$.
Then $P$ is the isotropy group of the flat metric $h\in \cM$ under the 
$\CDiff$ action.

If $(\varphi,\Omega)\in \CDiff$, \eqref{initial} shows that there is a
unique $p\in P$ such that to second order we have  
$\varphi = \varphi_p$ and $\Omega = \Omega_p$.
This defines a homomorphism $\CDiff \rightarrow P$ with kernel 
$\CDiff_0$.  Thus $\CDiff_0$ is a normal subgroup of $\CDiff$ and 
$\CDiff=P\cdot\CDiff_0$.  Moreover, $\cM/\CDiff_0$ has a natural 
left $P$-action.  

Just as in the pseudo-Riemannian case, local conformal invariants
correspond precisely to $P$-invariants of  
$\cM/\CDiff_0$.  The conformal jet isomorphism theorem provides a
tensorial description of $\cM/\CDiff_0$.  Since  
$P\subset GL(n+2,\R)$, $P$ acts on tensor powers of $\R^{n+2}$, not
$\R^n$.  Thus one anticipates a description as a $P$-space in terms of
tensors in $n+2$ dimensions.    

A significant difference from the pseudo-Riemannian case which will appear 
below is the fact 
that the structure of $\cM/\CDiff_0$ depends in a fundamental way on
whether $n$ is 
even or odd.  This is not evident at a superficial level.  The tangent
space $T(\cM/\CDiff_0)$ 
is isomorphic to the
quotient of a particular dual generalized Verma module, the jets of
trace-free symmetric 2-tensors of weight 2, by the image of the conformal
Killing operator acting on jets of vector fields.  (See Lemma~\ref{iso1}
below.)  These spaces have natural structures as 
$(\mathfrak{g},P)$-modules, where $\mathfrak{g}$ denotes the 
Lie algebra of $G$.  As a $(\mathfrak{g},P)$-module, 
this quotient is irreducible if $n$ is odd, but has a unique proper
$(\mathfrak{g},P)$-submodule with irreducible quotient if $n$ is even.   
Geometrically, the distinction is exhibited by the existence of the ambient  
obstruction tensor, a conformally invariant natural tensor which exists
only in even dimensions. 

Next we formulate the jet isomorphism theorem for conformal geometry for
$n$ odd.  
\begin{definition}\label{Rt}
Let $n\geq 3$ be odd.  
The space $\cRt\subset \prod_{r=0}^\infty
\bigwedge^2\R^{n+2}{}^*\otimes\bigwedge^2\R^{n+2}{}^*\otimes
\bigotimes^r\R^{n+2}{}^*$ 
is the set of lists $(\Rt^{(0)},\Rt^{(1)},\Rt^{(2)},\cdots )$ with 
$\Rt^{(r)}\in\bigwedge^2\R^{n+2}{}^*\otimes\bigwedge^2\R^{n+2}{}^*\otimes  
\bigotimes^r\R^{n+2}{}^*$, such that:   
\begin{enumerate}
\item
$\Rt_{I[JKL],M_1\cdots M_r}=0$
\item
$\Rt_{IJ[KL,M_1]M_2\cdots M_r} = 0$
\item
$\Rt_{IJKL,M_1\cdots [M_{s-1}M_s]\cdots M_r}= {\tilde
   Q}^{(s)}_{IJKLM_1\cdots  M_r}(\Rt)$    
\item
$\htt^{IK}\Rt_{IJKL,M_1\cdots M_r} = 0$
\item
$\Rt_{IJK0,M_1\cdots M_r}= 
-\sum_{s=1}^r \Rt_{IJKM_s,M_1\cdots \widehat{M_s} \cdots M_r}$.
\end{enumerate}
Here ${\tilde Q}^{(s)}_{IJKLM_1\cdots 
M_r}(\Rt)$ is the same quadratic expression in the 
$\Rt^{(r')}$, $r'\leq r-2$, as in Definition~\ref{R},  
except that now the tensors are the $\Rt^{(r')}$ instead of the $R^{(r')}$
and 
the contractions  are taken with respect to $\htt$ instead of $h$.
Condition (5) in case $r=0$ is interpreted as $\Rt_{IJK0}= 0$.    
\end{definition}

Conditions (1)--(4) are invariant under all of $O(\htt)$.  Condition (5) is
certainly not invariant under $O(\htt)$, but it is almost invariant under
$P$.  Recall that
$p\in P$ given by \eqref{pform} satisfies $pe_0 = ae_0$, $a>0$.  So (5)  
is invariant under $P$ modulo the rescaling of $e_0$.  To correct the
scaling, for 
$w\in \C$ define the character $\sigma_w:P\rightarrow \C$ by 
$\sigma_w(p)=a^{-w}$.  Then if we define the $P$ action by viewing 
\begin{equation}\label{Pstruc}
\cRt\subset \prod_{r=0}^\infty\,\,\textstyle{
\bigwedge^2\R^{n+2}{}^*\otimes\bigwedge^2\R^{n+2}{}^*\otimes
\bigotimes^r\R^{n+2}{}^*\otimes \sigma_{-2-r}}, 
\end{equation}
then the scaling of the $\Rt^{(r)}$ cancels the scaling of $e_0$ and  
condition (5) becomes invariant under $P$.  Thus $\cRt$ becomes a
$P$-space.  (In the factor $\sigma_{-2-r}$, the $-2$ could be replaced by
any other number and condition (5) would still be $P$-invariant.  The
choice of $-2$ is necessary for the map $c$ below to be
$P$-equivariant.)  The conformal jet isomorphism theorem for $n$ odd is
then the following.

\begin{theorem}\label{jetisoodd}
If $n$ is odd, then there is a $P$-equivariant polynomial bijection 
$c:\cM/\CDiff_0\rightarrow \cRt$ with polynomial inverse.
\end{theorem}

If $n$ is even, the analogous statement holds only for truncated jets:
there is a bijection from $(n-1)$-jets of 
metrics mod $\CDiff_0$ to a correspondingly truncated version of the space
$\cRt$.  An infinite order extension of this result for $n$ even will be 
discussed in \S\ref{even}.  

The jet isomorphism theorem reduces the study of conformal invariants to
the study of $P$-invariants of $\cRt$.  This is important because algebraic
tensorial operations can be utilized to construct and study conformal
invariants.  

Next we discuss the origin of the space $\cRt$ and the construction of the
map $c$.  As described above, the conformal geometry of the 
quadric $\cQ$ naturally arises from the metric
$\htt_{IJ}dx^Idx^J$ on $\R^{n+2}$.   In \cite{FG1}, a version of the metric 
$\htt$ for a general conformal manifold, called the ambient metric, was
introduced.  The tensors $\Rt^{(r)}$ arise as the iterated covariant
derivatives of the curvature tensor of the ambient metric.  

Suppose that $M$ is a smooth manifold with a conformal class $[g]$ of
metrics of signature $(p,q)$.  
The metric bundle $\cG$ of $[g]$ is 
$\cG=\{(x,t^2g(x)): x\in M,\,t>0\}\subset \bigodot^2T^*M$, where $g$ is a  
metric in the conformal class.  The fiber variable $t$ on $\cG$ is
associated to the 
metric $g$ and provides an identification $\cG\cong \R_+\times M$.   
There is a tautological symmetric $2$-tensor ${\bf g_0}$ on $\cG$ defined
by  
${\bf g_0}(Y,Z)=\underline{g}(\pi_*Y,\pi_*Z)$, where $\pi:\cG\rightarrow M$
is the projection 
and $Y$, $Z$ are tangent vectors to $\cG$ at $(x,\underline{g})\in \cG$. 
The family of dilations $\delta_s:\cG\rightarrow \cG$ defined by  
$\delta_s(x,\underline{g})= (x,s^2\underline{g})$ defines an $\R_+$ action
on $\cG$, and one has  
$\delta_s^*{\bf g_0} = s^2{\bf g_0}$.   
We denote by $T= \frac{d}{ds} \delta_s |_{s=1}$ the vector field on $\cG$
which is the infinitesimal generator of the dilations $\delta_s$.   Note
that ${\bf g_0}$ is degenerate:  ${\bf g}_0(T,Y)=0$ for all $Y\in T\cG$.   
In the case that $(M,[g])$ is the quadric $\cQ$ with its conformal
structure defined above, $\cG$ can be identified with $\cN/\{\pm I\}$, 
${\bf g}_0$ with $\htt|_{T\cN}$, and $T$ with $X$.  

The ambient space is defined to be $\cGt= \cG\times \R$;  
the coordinate in the $\R$ factor is typically written $\rho$.
The dilations $\delta_s$ extend to $\cGt$ acting on the $\cG$ factor
and we denote also by $T$ the infinitesimal generator of the $\delta_s$ on 
$\cGt$.  
We embed $\cG$ into $\cGt$ by $\iota: z\rightarrow (z,0)$ for $z\in
\cG$, and we identify $\cG$ with its image under $\iota$.  
We say that a subset of $\cGt$ is homogeneous if it is invariant under
the $\delta_s$ for all $s>0$.  We say that a map between 
homogeneous subsets of $\cGt$ is homogeneous if it commutes with the
$\delta_s$.   

\begin{definition}\label{gt}
Let $n$ be odd.  An ambient metric $\gt$ for $(M,[g])$ is a smooth 
metric of signature $(p+1,q+1)$ on a homogeneous
neighborhood of $\cG$ in $\cG\times \R$ satisfying:
\begin{enumerate}
\item
$\delta_s^*\gt = s^2 \gt$
\item
$\iota^* \gt = {\bf g_0}$
\item
$\Ric(\gt)=0$ to infinite order along $\cG$.
\end{enumerate}
\end{definition}
\noindent
The main result concerning existence and uniqueness of the ambient metric
for $n$ odd is:
\begin{theorem}
If $n$ is odd, there exists an ambient metric for $(M,[g])$.  It is unique
up to:
\begin{enumerate}
\item[(a)] Pullback by a homogeneous diffeomorphism $\Phi$ satisfying
$\Phi|_\cG = I$, and 
\item[(b)]
Addition of a tensor homogeneous of degree 2 which vanishes to infinite
order along $\cG$.         
\end{enumerate}
\end{theorem}

The proof proceeds by the introduction of a gauge normalization to break
the diffeomorphism invariance together with a formal power series   
analysis of the equations $\Ric(\gt)=0$.  See \cite{FG2}. 

When $n\geq 4$ is even, there is an obstruction at order $n/2$ to existence
of a metric satisfying (1)--(3), which is a conformally invariant natural
tensor called the ambient 
obstruction tensor.  However, there is a solution up to this order, again
unique up to homogeneous diffeomorphism and up to a term vanishing to 
order $n/2$.    

The solution has an extra geometric property:  for each $p\in \cGt$, the
parametrized dilation orbit $s\rightarrow \delta_sp$ is a geodesic for
$\gt$ (to infinite order along $\cG$).  

The diffeomorphism ambiguity in $\gt$ can be fixed by the choice of a
metric $g$ in the conformal class.  As described above, the choice of such
a metric determines an identification $\cG\cong \R_+\times M$, and
therefore an identification $\cGt\cong \R_+\times M \times \R$.  
\begin{definition}\label{normalform}
A metric $\gt$ satisfying conditions (1) and (2) in Definition~\ref{gt} is
said to be in normal form relative to $g$ if in  
the identification $\cGt\cong \R_+\times M \times \R$ induced by $g$, one
has 
\begin{enumerate}
\item
$\gt = 2t\,dt\,d\rho + {\bf g_0}$ at $\rho = 0$, and
\item
The curve
$\rho\rightarrow (t,x,\rho)$ is a geodesic for $\gt$ for each choice of  
$(t,x)\in \R_+\times M$.
\end{enumerate}
\end{definition}
\noindent
If $n$ is odd, an ambient metric $\gt$ can always be found which is in  
normal form relative to $g$, and it is uniquely determined up to 
$O(\rho^{\infty})$.  Each term in the Taylor expansion at $\rho =0$ of such
a $\gt$ in normal form relative to $g$ is given by a polynomial expression
in $g^{-1}$ and in derivatives of $g$.  

An analogue in conformal geometry of the curvature tensor 
and its covariant derivatives for a pseudo-Riemannian metric  
are the restrictions to $\cG$ of the
curvature tensor and its covariant derivatives of the ambient metric.  
These can be 
interpreted as sections of weighted tensor powers of the cotractor bundle  
associated to the conformal structure; see \cite{CG}, \cite{BG} and
\cite{FG2}.  For our purposes,  
the map $c$ in Theorem~\ref{jetisoodd} can be defined 
directly as follows.  For $g\in \cM$,   
choose a metric also denoted $g$ defined near $0\in \R^n$
with the prescribed Taylor expansion.  There is an ambient   
metric $\gt$ in normal form relative to $g$,  
uniquely determined to infinite order in $\rho$.  Define the tensors
$\Rt^{(r)}$ in Theorem~\ref{jetisoodd} to be the iterated covariant
derivatives of the curvature tensor of $\gt$ evaluated   
at $t=1$, $x=0$ and $\rho = 0$.  It can be shown that these covariant
derivatives satisfy the relations (1)--(5) in Definition~\ref{Rt} which
define $\cRt$.  
Relations (1)--(3) hold for the covariant derivatives of curvature of any
metric.  Relation (4) follows from the fact that $\gt$ is Ricci-flat to
infinite order, and (5) is a consequence of the homogeneity of
$\gt$ and the fact that the dilation orbits are geodesics to infinite
order.  
Now using the fact that the ambient curvature tensors are 
tensors on the ambient space, it can be shown that this map
$\cM\rightarrow \cRt$ passes to a 
map $c:\cM/\CDiff_0\rightarrow \cRt$ which is $P$-equivariant.   
Details can be found in \cite{FG2}.  

The invertibility of $c$ in Theorem~\ref{jetisoodd} is also proved in
\cite{FG2}.  As in the pseudo-Riemannian case, one first constructs a   
slice for the $\CDiff_0$ action, using geodesic normal coordinates and a
``conformal normal form'' which normalizes away the freedom of the 
derivatives of the conformal factor of order two or more.  (Actually, the
formulation of the jet isomorphism theorem in \cite{FG2} is in terms of
this slice rather than in terms of the space $\cM/\CDiff_0$.)  A
linearization argument reduces the matter to showing the invertibility of
the linearization $dc$ of $c$ at the flat metric $h$.      
Then the main part of the proof consists of an algebraic study 
of the relations obtained by linearizing (1)--(5) (they are all already
linear except for (4)) and a       
direct analysis of $dc$.  A more conceptual
proof of the invertibility of the linearized map will be outlined   
in the next section as an application of the results on the ambient
realization of the deformation complex.  

\section{Ambient lift of deformation complex}\label{deform}

In this section we introduce the conformal deformation complex and indicate 
how it may be realized ambiently in odd   
dimensions.  We then sketch a proof of the invertibility of the map $c$ in
Theorem~\ref{jetisoodd} using the exactness of the deformation 
complex on jets together with this ambient realization.  Details will
appear in \cite{GH2}.            

Recall from the previous section that the conformal group $G=O(\htt)/\{\pm
I\}$ acts conformally on the quadric   
$\cQ$ with isotropy group $P$ so that $\cQ=G/P$, and that there is an
embedding $i:\R^n\to \cQ$ as an open dense subset on which the
conformal structure is represented by the flat metric $h$.  To each
finite-dimensional representation of $P$ is associated a homogeneous vector
bundle on $\cQ=G/P$ and therefore also on $\R^n\hookrightarrow\cQ$.  
Familiar examples include:
\medskip
\begin{itemize}
\item
$\cD_w$, $w\in \C$:  the bundle of conformal densities of weight $w$,  
induced by $\sigma_w$
\item
$T\cQ$:  the tangent bundle, induced by $p\mapsto a^{-1}m$
\item
$\bigwedge^r$, $0\leq r\leq n$:  the bundle of $r$-forms (the
$r$-th exterior power of the cotangent bundle).    
\end{itemize}
\medskip
Set $\bigwedge^r(w)=\bigwedge^r\otimes \cD_w$.  We denote by $\cE(w)$,
$\cE^r(w)$ the sheaf of germs of smooth sections of $\cD_w$, 
$\bigwedge^r(w)$, resp., and by $\cE_\cU(w)$, $\cE^r_\cU(w)$ the space of
sections on an open set $\cU\subset \cQ$.  

For $0\leq s \leq r\leq n$, define $\bigwedge^{r,s}$ 
to be the homogeneous bundle on $\cQ$ of   
covariant tensors with Young symmetry given by the Young diagram  
\begin{equation}\label{young}
\begin{array}{c}
\begin{picture}(30,70)(-5,0)
\put(-25,33){\mbox{$r \left\{
\begin{picture}(0,38)\end{picture}
\right.$}}
\put(25,48){\mbox{$\left.
\begin{picture}(0,23)\end{picture}
\right\}s$}}
\put(0,0){\line(0,1){70}}
\put(10,0){\line(0,1){70}}
\put(20,30){\line(0,1){40}}
\put(0,70){\line(1,0){20}}
\put(0,60){\line(1,0){20}}
\put(3,45){\mbox{$\vdots$}}
\put(13,45){\mbox{$\vdots$}}
\put(0,40){\line(1,0){20}}
\put(0,30){\line(1,0){20}}
\put(3,15){\mbox{$\vdots$}}
\put(0,10){\line(1,0){10}}
\put(0,0){\line(1,0){10}}
\end{picture}
\end{array}
\end{equation}
Explicitly, $\bigwedge^{r,s}$ is the subbundle of 
$\bigwedge^r\otimes\bigwedge^s$ consisting of those tensors
$$
f_{i_1\dots i_rj_1\dots j_s}=f_{[i_1\dots i_r][j_1\dots j_s]}\in 
\textstyle{\bigwedge^r\otimes\bigwedge^s}
$$
which satisfy
$$
f_{[i_1\dots i_rj_1]j_2\dots j_s}=0.  
$$
Note that $\bigwedge^{r,0}=\bigwedge^r$ and that 
$\bigwedge^{1,1}=\bigodot^2$ is the bundle of symmetric
2-tensors.   
We denote by $\bigwedge^{r,s}_0\subset \bigwedge^{r,s}$ the subbundle of
tensors which 
are trace-free with respect to a metric in the conformal class, by 
$\bigwedge^{r,s}(w)$, $\bigwedge^{r,s}_0(w)$ the respective tensor products
with $\cD_w$, and by $\cE^{r,s}(w)$, $\cE^{r,s}_0(w)$ the sheaves of germs
of sections.  Each of the bundles $\bigwedge^{r,s}_0(w)$ is an irreducible  
homogeneous bundle; i.e., it is induced by an irreducible representation of
$P$.  

We will represent sections of $\bigwedge^{r,s}(w)$ in either of two ways.   
On $\R^n\hookrightarrow \cQ$, we can use $h$ to trivialize the  
density bundle and can thereby identify a section with a tensor  
field $u$ on an open subset of $\R^n$ having the symmetries indicated
above.  Alternately, we can view a section 
as a homogeneous tensor field $f$ on an open subset of $\cN$.  Define the 
dilations 
$\delta_\lambda:\R^{n+2}\to \R^{n+2}$ by 
$\delta_\lambda(x)=\lambda x$ for $\lambda\in \R\setminus\{0\}$.  
Then for $\cU\subset \cQ$ open, there is a 1-1 correspondence between
$\cE^{r,s}_\cU(w)$ and the set of smooth sections $f$ of
$\bigotimes^{r+s}T^*\cN$ on $\pi^{-1}(\cU)$
which have the symmetries above and which satisfy  
\begin{equation}\label{2eqns}
\delta_{\lambda}^*f = |\lambda|^wf,\qquad X\into f =0.
\end{equation}
Here the condition $X\into f=0$ is interpreted to mean that the contraction
of $X$ into  
every index of $f$ vanishes.  

We now work on $\R^n$, viewed as a subset of $\cQ$, and use its usual  
coordinates and the flat metric $h$.  Define differential operators  
\begin{align}\label{diffops}
\begin{aligned}
&d_1:\cE^{r,s}\rightarrow \cE^{r+1,s}\\ 
&d_2:\cE^{r,s}\rightarrow  
\cE\left(\textstyle{\bigwedge}^r\otimes \textstyle{\bigwedge}^{s+1}\right)
\\  
&\delta_1:\cE^{r,s}\rightarrow 
\cE\left(\textstyle{\bigwedge}^{r-1}\otimes \textstyle{\bigwedge}^s\right)
\\  
&\delta_2:\cE^{r,s}\rightarrow \cE^{r,s-1} 
\end{aligned}
\end{align}
by:
\begin{align*}
\begin{aligned}
(d_1 u)_{i_0i_1\cdots i_r j_1\cdots j_{s}}&
=\partial_{[i_0} u_{i_1\cdots i_r] j_1\cdots j_{s}}\\
(d_2 u)_{i_1\cdots i_r j_0\cdots j_s}&
=\partial_{[j_0} u_{|i_1\cdots i_r |j_1\cdots j_{s}]}\\
(\delta_1 u)_{i_1\cdots i_{r-1} j_1\cdots j_s}&
=-\partial^{k} u_{i_1\cdots i_{r-1}k j_1\cdots j_s}\\
(\delta_2 u)_{i_1\cdots i_r j_1\cdots j_{s-1}}&
=-\partial^{k} u_{i_1\cdots i_r j_1\cdots j_{s-1}k}.    
\end{aligned}
\end{align*}
Here the $|i_1\cdots i_r|$ indicates indices excluded from the  
skew-symmetrization.  The derivatives are coordinate derivatives on $\R^n$ 
and the contractions are with respect to $h$.  In making this definition,
we momentarily ignore   
the weights and the structure as homogeneous bundles and view these simply
as differential operators on tensor fields.      

For $n\geq 4$, the deformation complex on $\R^n$ is:
\begin{equation}\label{defcx}
\begin{split}
  0\to\mathfrak{g}\to\cE^{1}(2)&\stackrel{D_{0}}{\longrightarrow} 
  \cE_0^{1,1}(2)\stackrel{D_{1}}{\longrightarrow}
\cE_0^{2,2}(2)\stackrel{D_{2}}{\longrightarrow}
\cE_0^{3,2}(2)\\
&\to\cdots\to\cE_0^{n-2,2}(2)\stackrel{D_{n-2}}{\longrightarrow}
  \cE_0^{n-1,1}\stackrel{D_{n-1}}{\longrightarrow}
\cE^{n-1}(-2)\to 0, 
\end{split}
\end{equation}
where 
$$
\begin{aligned}
D_0&=\tf\text{Sym}\,d_2\\
D_1&=\tf d_1d_2
\\
D_r&= \tf d_1
\quad r=2,3,\dots, n-3\\ 
D_{n-2}&
=\delta_2d_1 
\\
D_{n-1}&=\delta_2. 
\end{aligned}
$$
\noindent
Here $\tf$ denotes the trace-free part with respect to $h$ and 
$\text{Sym}$ denotes symmetrization over the two indices.  
When $n=4$, the $\cE_0^{2,2}(2)$ on the first line and the 
$\cE_0^{n-2,2}(2)$ on the second line are the same space, so 
$D_2=\delta_2d_1$ maps into $\cE_0^{3,1}$ and the space 
$\cE_0^{3,2}(2)$ does not occur.  In higher dimensions, the spaces
between $\cE_0^{2,2}(2)$ and $\cE_0^{n-1,1}$ are the
$\cE_0^{r,2}(2)$ for $3\leq r\leq n-2$.  The $D_r$ are the expressions on
$\R^n$ of 
$G$-equivariant differential operators between the indicated homogeneous
vector bundles on $G/P$, or equivalently between the sheaves of their germs
of local 
sections.  The space $\mathfrak{g}$ is the locally constant sheaf.  
The bundle $\bigwedge^1(2)$ is isomorphic to the tangent bundle by raising 
the index, and in this realization the map $\mathfrak{g}\to \cE^1(2)$
is the infinitesimal $G$-action.   

The deformation complex is a complex, i.e. the composition of two
successive operators vanishes.  It can be thought of as analogous to the
deRham complex; it has the same length as  the deRham complex.  The 
operators $D_1$ and $D_{n-2}$ are second order; all other $D_r$ are first
order.  The deformation complex was constructed explicitly ``by hand'' by
Gasqui-Goldschmidt in \cite{GG} on a general conformally flat manifold.  In
the homogeneous case it is a particular case of a generalized  
Bernstein-Gelfand-Gelfand complex (see \cite{L} for the introduction of
gBGG complexes in the algebraic setting).  In the 3-dimensional case 
the deformation complex takes a special form: 
\begin{equation}\label{dim3}
  0\to\mathfrak{g}\to\cE^{1}(2)\stackrel{D_{0}}{\longrightarrow} 
  \cE_0^{1,1}(2)\stackrel{D_{1}}{\longrightarrow}
\cE_0^{2,1}\stackrel{D_{2}}{\longrightarrow}
\cE^{2}(-2)\to 0,
\end{equation}
where $D_0$ is as above, $D_2=\delta_2$, and 
$D_1=\tf \delta_2d_1d_2 $ is third order.   

The main fact that we will need about the deformation
complex is that it is exact on jets; i.e., if an infinite-order jet of a 
section of one of 
the bundles at a point is annihilated by the corresponding 
operator as a jet, then it is in the image of the previous 
operator acting on jets at that point.  This fact is proved in \cite{GG} 
and is also contained in the theory of the generalized BGG complexes.     

This complex is called the deformation complex because its first terms
describe the infinitesimal deformation of conformal structures.   
The first operator $D_0$ corresponds to the conformal Killing operator  
$\tf\mathcal{L}_Vh$, where $\mathcal{L}$ denotes the 
Lie derivative and $V$ a vector field, which is obtained by linearizing the
action of 
diffeomorphisms on conformal structures.  Its kernel $\mathfrak{g}$
consists of the infinitesimal conformal transformations.  For $n\geq 4$,
the operator 
$D_1$ is the linearization of the map which takes the Weyl tensor of a
metric, and $D_2$ is the linearization of the Bianchi identity satisfied by
such a Weyl tensor.  For $n=3$, $D_1$ is the linearization of the Cotton
tensor, and $D_2$ the linearization of the ``Bianchi identity'' satisfied
by a Cotton tensor of a metric.

We wish to give an alternate description of the deformation complex for 
$n$ odd in which the spaces and maps are defined on the  ambient space.  
Other descriptions and curved versions are contained in \cite{CSS} and 
\cite{CD} (in much greater generality), and in \cite{GP}.     
We begin by introducing the ambient versions of the spaces appearing in the  
complex.  

For $0\leq s\leq r$, denote by $\wt{}^{r,s}$ the vector bundle of tensors
of rank $r+s$ on 
$\R^{n+2}$ having the Young symmetry \eqref{young} and by $\wt{}^{r,s}_0$  
the subbundle of those tensors which are trace-free with respect to
$\htt$.  We write 
$\tilde{d}_1$, $\tilde{d}_2$, $\tilde{\delta}_1$, $\tilde{\delta}_2$
for the operators on $\R^{n+2}$ analogous to \eqref{diffops} and 
$\tilde{\Delta}=\htt^{IJ}\pa^2_{IJ}$ for the Laplacian with
respect to $\htt$.  $\tilde{\Delta}$ acts on sections of 
$\wt{}^{r,s}$ componentwise with respect to the standard basis.   
Recall that $X=x^I\pa_I$ denotes the Euler field on 
$\R^{n+2}$, whose components are thus given by $X^I=x^I$.  

Let $\pi:\R^{n+2}\setminus \{0\}\to
\mathbb{P}^{n+1}$ be the projection.
Let $0\leq s\leq r$ and $w\in \C$.  For $\cV\subset \mathbb{P}^{n+1}$ open, 
define 
$\cEt{}^{r,s}_{\cV}(w)$ to be the space of sections $\ft$ of $\wt{}^{r,s}$   
on $\pi^{-1}(\cV)$ which satisfy  
$\delta_\lambda^*\ft = |\lambda|^w\ft$ for $\lambda\in \R\setminus\{0\}$, 
and $\cEt{}^{r,s}_{0,\cV}(w)$ to be the subspace of trace-free sections.   
The assignments $\cV\to \cEt{}^{r,s}_{\cV}(w)$, $\cEt{}^{r,s}_{0,\cV}(w)$  
define presheaves on $\mathbb{P}^{n+1}$ whose associated sheaves we denote
by $\cEt{}^{r,s}(w)$, $\cEt{}^{r,s}_0(w)$, resp.  
Observe that pullback defines a natural action of $O(\htt)$ on the total
space of these
sheaves and $\pm I$ acts by the identity, so that
$G=O(\htt)/\{\pm I\}$ also acts.  

Recall that $\cN$ is the null cone of $\htt$ and that 
$\pi:\cN\to\cQ$.  
If $\cU\subset \cQ$ is open, define
$\cH{}^{r,s}_{\cU}(w)$ to be the space of infinite-order jets along 
$\pi^{-1}(\cU)$  
of sections $\ft\in \cEt{}^{r,s}_{0,\cV}(w)$ for some $\cV\subset 
\mathbb{P}^{n+1}$ open, $\cU\subset \cV$,  
which satisfy the following  
equations formally to infinite order along $\pi^{-1}(\cU)$: 
\begin{equation}\label{3eqns}
\tilde{\Delta}\ft = 0,\qquad \tilde{\delta}_1\ft = 0, \qquad X\into \ft=0.      
\end{equation}
Here again $X\into \ft=0$ is interpreted to mean that the contraction of
$X$ into any index of $\ft$ vanishes.  By the symmetries of $\ft$, this is
equivalent to $X^{I_1}\ft_{I_1\cdots I_rJ_1\cdots J_s}=0$.  Similarly, 
$\dt_1\ft=0$ implies $\dt_2\ft = 0$.  
The assignment $\cU\to \cH{}^{r,s}_{\cU}(w)$ defines a presheaf on $\cQ$,
whose associated sheaf we denote by $\cH{}^{r,s}(w)$.  Since the equations
\eqref{3eqns} are invariant under $O(\htt)$, the conformal group $G$ acts
on $\cH{}^{r,s}(w)$ covering the translations of
$G$ on $\cQ=G/P$.  Thus $\cH{}^{r,s}(w)$ is a ``homogeneous sheaf'' on
$G/P$ in the 
same sense as in the definition of a homogeneous vector bundle.    

The ambient realization for $\cE^{r,s}_0(w)$ in the unobstructed cases
is given by the following theorem.  
\begin{theorem}\label{amlift}
Suppose $n\geq 3$.  Let $0\leq s\leq r\leq n$ and $w\in \C$.    
\begin{itemize}
\item
If $r>s=0$, assume that $w\neq 2r-n$.
\item
If $s>0$, assume that $w\neq r+2s-n-1$, $2r+s-n$.
\end{itemize}
If $w+n/2-r-s\notin \N$, then 
$$
\cE^{r,s}_0(w)\cong \cH{}^{r,s}(w)
$$
$G$-equivariantly as sheaves on $\cQ$.  
\end{theorem}

The disallowed values correspond to the existence of certain particular
$G$-invariant 
differential operators which act on $\cE^{r,s}_0(w)$.  In particular, 
if the dual generalized Verma module associated to $\cE^{r,w}_0(w)$ is
irreducible as a $(\mathfrak{g},P)$-module, then Theorem~\ref{amlift}
applies to $\cE^{r,s}_0(w)$.  It is important to note, 
as we will see, that not all $G$-invariant differential operators obstruct
the isomorphism asserted by Theorem~\ref{amlift}.  Otherwise stated,
Theorem~\ref{amlift} applies to many homogeneous bundles $\cE^{r,s}_0(w)$
for which the associated dual Verma modules are not irreducible as 
$(\mathfrak{g},P)$-modules.  

The map in one direction in Theorem~\ref{amlift} is evident, and exists for 
all values of $w$, $r$ and $s$.  For $\cU\subset \cQ$ open,
view elements of $\cE^{r,s}_{0,\cU}(w)$ as    
covariant tensor fields $f$ on $\pi^{-1}(\cU)\subset \cN$ satisfying  
\eqref{2eqns} as described above.  Let
$\iota:\cN\to\R^{n+2}$ denote the inclusion.  If 
$\ft\in \cEt{}^{r,s}_{0,\cV}(w)$ for some $\cV\subset 
\mathbb{P}^{n+1}$ open with $\cU\subset \cV$, and $\ft$  
satisfies $(X\into \ft)|_{\pi^{-1}(\cU)} =0$, it is clear that    
$f=\iota^*\ft$ satisfies \eqref{2eqns}.  One checks easily that  
$f$ also satisfies the trace-free condition with respect to $h$ so that   
$f$ defines an element of $\cE^{r,s}_{0,\cU}(w)$.  Passing to jets along 
$\pi^{-1}(\cU)$  
and restricting consideration to $\cH^{r,s}(w)$ gives a $G$-equivariant
map   
$$
\iota^*:\cH^{r,s}(w)\rightarrow \cE_0^{r,s}(w).  
$$
The content of Theorem~\ref{amlift} is that under the stated restrictions
on the parameters, this map is a isomorphism.  That is, 
each element of $\cE^{r,s}_{0,\cU}(w)$ has a unique extension (ambient
lift) as an element of $\cH^{r,s}_{\cU}(w)$.    

There are two main steps in the proof of Theorem~\ref{amlift}.  The first 
is called the ``initial lift'', and corresponds to defining on
$\pi^{-1}(\cU)$ the 
components of $\ft$ transverse to $\cN$ to obtain a section of   
$\Lt{}^{r,s}_0|_{\pi^{-1}(\cU)}$ homogeneous of degree $w$.  The ideas in
this step  
go back to Tracy Thomas for special cases of the symmetries including
differential forms; he called the process of defining the transverse   
components ``completeing'' the tensor.  This step is closely related to
what are nowadays called differential splittings, about which there
is a substantial literature.  
This first step leads to the excluded values of the parameters indicated in
the bullets above.  
The second step involves ``harmonic'' extension of the completed tensor to
higher order off $\cN$ in such a way as to make the equations
\eqref{3eqns} as well as the trace-free condition hold to infinite order.  
The condition $w+n/2-r-s\notin \N$ arises in this second step.  Part of the
difficulty of the proof, especially for more complicated symmetries, is
making sure that the steps can be carried out 
consistently so that all the required conditions hold to all orders.  

To give an idea how this works, we sketch the details for the scalar case 
$r=s=0$ 
and the case $r=1$, $s=0$ of 1-forms.  The scalar case is studied in detail
in \cite{EG}.  I am grateful to
M. Eastwood for providing the argument below in the case of 1-forms.   

In the case $r=s=0$, Theorem~\ref{amlift} asserts that $\mathcal{E}(w)\cong 
\cH(w)$ if $w+n/2 \notin \N$, where $\cE(w)$ denotes the sheaf of germs of 
densities of weight $w$ on $\cQ$ and $\cH(w)$ denotes the sheaf of jets
along $\cN$ 
of homogeneous functions of degree $w$ which satisfy $\tilde{\Delta}\ft=0$
to infinite order.  Set $Q=\htt_{IJ}x^Ix^J$; then this is the same as 
showing that given $f$ homogeneous of degree $w$ on
$\pi^{-1}(\cU)\subset\cN$, there exists 
a unique infinite order jet $\ft$ homogeneous of degree $w$ satisfying
$\tilde{\Delta}\ft = O(Q^\infty)$ and $\ft|_{\pi^{-1}(\cU)}=f$.  The 
initial lift 
step is vacuous in this case.  For the harmonic extension step, the Taylor  
expansion of $\ft$ is constructed inductively.  A key observation is that 
\begin{equation}\label{comm}
[\tilde{\Delta},Q^k]=2kQ^{k-1}(2X+n+2k). 
\end{equation}
Suppose that $\ft^{(k)}$ has been constructed which satisfies
$\Dt\ft^{(k)}=O(Q^{k-1})$.  Set 
$$
\ft^{(k+1)}=\ft^{(k)}+Q^k\eta \qquad \text{for}\quad \eta\in \cEt(w-2k).   
$$
Then 
\[
\begin{split}
\Dt\ft^{(k+1)} &= \Dt\ft^{(k)} +\Dt(Q^k\eta)\\
&=\Dt\ft^{(k)} + [\Dt,Q^k]\eta +O(Q^k)\\  
&=\Dt\ft^{(k)} + 2kQ^{k-1}(2X+n+2k)\eta + O(Q^k)\\ 
&=\Dt\ft^{(k)} + 2k(n+2w-2k)\eta Q^{k-1}+ O(Q^k).  
\end{split}
\]
So if $n+2w\neq 2k$, $\eta$ can be uniquely chosen so that 
$\Dt\ft^{(k+1)}=O(Q^{k})$.  Thus if $w+n/2\notin\N$, then   
the induction can be carried out to all orders.  

If $n/2+w =m\in \mathbb{N}$, then harmonic extension is obstructed by 
the conformally invariant operator $\Delta^m=(h^{ij}\pa^2_{ij})^m$ on  
$\R^{n}$.  

Consider now the case $r=1$, $s=0$.  Theorem~\ref{amlift} asserts that if
$w\neq 2-n$ and $w+n/2-1\notin \N$, then 
$\cE^1(w)\cong \cH{}^1(w)$, where we have written $\cH{}^1(w)$ for 
$\cH{}^{1,0}(w)$.  Recall that $\ft\in \cH{}^1(w)$ means that 
$\ft$ is a jet of a section of $\cEt{}^1(w)$ satisfying   
the equations \eqref{3eqns} to infinite order.
We write $\dt$ for $\dt_1$ since $\dt_2$ vanishes in this case.     

Given a 1-form $f$ on $\pi^{-1}(\cU)$ which is homogeneous of degree $w$
and 
satisfies $f(X) =0$, we can choose some $\ft$ which is a section in 
$\cEt{}^1_\cV(w)$ for some $\cV\supset \cU$ such that $\iota^*\ft =f$.
Such an $\ft$ is uniquely 
determined up to addition of $\psi dQ +Q\phi$ with $\psi$ a function and 
$\phi$ a 1-form, both of homogeneity $w-2$.  We can certainly choose $\ft$
to start with so that $\ft(X) = O(Q^2)$; in fact we could make  
$\ft(X) = O(Q^\infty)$, but $O(Q^2)$ will suffice.  Now try to 
determine $\psi$, $\phi$ to maintain this condition on vanishing of
$\ft(X)$:   
$$
(\ft + \psi dQ +Q\phi)(X)=O(Q^2)
$$
and also to make
$$
\dt(\ft + \psi dQ +Q\phi)=O(Q).  
$$
The first equation gives $\quad\psi dQ(X) + Q\phi(X) = O(Q^2),\quad$ so 
$$
2\psi + \phi(X)=O(Q).
$$
The second equation gives 
$$
\dt\ft -2(n+w)\psi -2\phi(X) =O(Q),
$$
so combining gives
$$
\dt\ft -2(n+w-2)\psi  =O(Q). 
$$
If $n+w\neq 2$, this uniquely determines $\psi \mod Q$.  Then 
$(\ft+\psi dQ)|_{\pi^{-1}(\cU)}$ is the initial lift.      
Rename $\ft+\psi dQ$ to be a new $\ft$.    

Now all components of $\ft|_{\pi^{-1}(\cU)}$ have been determined.   
Write $\ft = \ft_Idx^I$; then each $\ft_I$ is a scalar function homogeneous  
of degree 
$w-1$.  Since $w+n/2-1\notin \N$, by the scalar case we can uniquely 
extend each $\ft_I$ harmonically to 
infinite order, and this is equivalent to the condition that $\ft$ 
satisfies $\Dt\ft=0$ to infinite order.  In particular, we conclude the 
uniqueness of an extension satisfying \eqref{3eqns} to infinite order.  For
existence, we  
claim that this harmonic extension automatically satisfies 
$\dt \ft=0$ and $\ft(X)=0$ to infinite order.  One first checks  
that the harmonic extension recovers the conditions $\dt \ft =O(Q)$ and  
$\ft(X) = O(Q^2)$ imposed above.  Then $\Dt$ and $\dt$ commute since 
they are constant coefficient operators on $\R^{n+2}$, so  
$\Dt \dt \ft=0$ to infinite order.  But $\dt \ft$ has homogeneity $w-2$ and  
$w-2+n/2\notin \mathbb{N}$, so uniqueness for the scalar case implies that 
$\dt \ft =0$ to infinite order.  The argument that $\ft(X)=O(Q^\infty)$ 
is similar.  One has $\Dt(\ft(X))=O(Q^\infty)$ since $\Dt 
\ft=O(Q^\infty)$ and $\dt\ft=O(Q^\infty)$.  Now $\ft(X)$ is homogeneous 
of degree $w$.  Since $w+n/2-1\notin\N$, we can apply the usual statement
of uniqueness in the scalar case unless $w+n/2=1$.  If
$w+n/2=1$, the argument in the scalar case proves uniqueness for
densities which are $O(Q^2)$.  Thus $\ft(X)=O(Q^\infty)$ holds in 
general.   

For general $r$, $s$, the algebra of the initial lift and the
consistency verification is more complicated, but the basic
idea is the same.  When $r>s=0$, the operator $\delta_1$ is conformally
invariant for $w=2r-n$ and obstructs the initial lift.
If $s>0$, there are two invariant operators obstructing the initial lift,
giving rise to the two excluded values of $w$.  For $r>s>0$ the invariant  
operators are $\delta_2$, $\pi^{r-1,s}\delta_1$ for $w=r+2s-n-1$,
$2r+s-n$,  respectively, where $\pi^{r-1,s}$ is the Young projector onto    
$\bigwedge^{r-1,s}$.  If $r=s>0$, the obstructing invariant operators are  
$\delta_2$ for $w=3r-n-1$ and $\delta_1\delta_2$, an iterated divergence, 
for $w=3r-n$.  

Theorem~\ref{amlift} implies 
a corresponding isomorphism obtained by taking jets at a point. 
Define $\cJ^{r,s}(w)$ to be the space of infinite-order jets at $[e_0]\in
\cQ$ of sections of $\bigwedge^{r,s}(w)$, and 
$\cJ^{r,s}_0(w)$ to be the subspace of jets which are trace-free to
infinite order.  The $G$-action on 
the sheaf $\cE^{r,s}_0(w)$ induces a $(\mathfrak{g},P)$-module structure on  
$\cJ^{r,s}_0(w)$ dual to a generalized Verma module.  Define 
$\cJt{}^{r,s}(w)$ to be the space of infinite order jets at $e_0\in
\R^{n+2}$ of  
sections of $\wt{}^{r,s}$ which are homogeneous of degree $w$, and by 
$\cJt{}^{r,s}_0(w)$ the subspace of jets which are trace-free to infinite
order.  The $G$-action on 
$\cEt{}^{r,s}(w)$ induces $(\mathfrak{g},P)$-module
structures on $\cJt{}^{r,s}(w)$, $\cJt{}^{r,s}_0(w)$.  Define 
$\cJt{}^{r,s}_\cH(w)\subset\cJt{}^{r,s}_0(w)$ to be the submodule
consisting of those jets for which the equations 
\eqref{3eqns} hold to infinite order at $e_0$.    It follows from
Theorem~\ref{amlift} that if $r$, $s$, $w$ satisfy the  
restrictions of Theorem~\ref{amlift}, then 
$$
\cJ^{r,s}_0(w)\cong \cJt{}^{r,s}_\cH(w)  
$$ 
as $(\mathfrak{g},P)$-modules.   This statement can be regarded as a  
``jet isomorphism theorem for $\cJ^{r,s}_0(w)$''  
providing an ambient description of the dual generalized
Verma modules.  Upon expressing a jet in $\cJt{}^{r,s}_\cH(w)$ as the list   
of tensors which are the successive derivatives of the section, one can 
realize the $P$-action in the ambient description 
in terms of reweighted tensor representations analogous to \eqref{Pstruc}.
See \cite{EG} for further discussion.

Observe that $r$, $s$, $w\in \Z$ for all of the spaces $\cE^{r,s}_0(w)$
which occur in the deformation complex \eqref{defcx}, \eqref{dim3}.
Therefore, if $n$ is odd, then the condition 
$w+n/2-r-s\notin \N$ in Theorem~\ref{amlift} is automatic for these spaces.
One verifies easily that for $n$ odd, the bulleted conditions in 
Theorem~\ref{amlift} hold for all spaces which occur in the deformation
complex except for the next to last one, $\cE^{n-1,1}_0$, for which the 
second bulleted condition is violated.  This corresponds to the fact that 
the operator in the complex acting on this space is $\delta_2$, which
is precisely the operator 
obstructing the ambient lift on this space.  Even though invariant
operators act on 
the other spaces in the deformation complex, namely the operators occuring
in the complex itself, the only one obstructing ambient lifts is the one  
acting on $\cE^{n-1,1}_0$.  So Theorem~\ref{amlift} provides an ambient
description of all of the other spaces in the complex.  It is not difficult
to identify the differential operators on $\R^{n+2}$ which correspond
in this realization to the operators in the deformation complex.  One thus 
obtains:
\begin{theorem}\label{deflift}
Let $n$ be odd.  The deformation complex 
with last two spaces removed can be realized as:  
\def\tos{\!\!\to\!\!}
\def\ltos{\!\!\longrightarrow\!\!}
$$
\begin{array}{ccccccccccccc}
0&\tos&\mathfrak{g}&\tos&
\cH^{1}(2)&\stackrel{{\tilde{D}}_{0}}{\ltos} &  
  \cH^{1,1}(2)&\stackrel{{\tilde{D}}_{1}}{\ltos}&
\cH^{2,2}(2)&\stackrel{{\tilde{D}}_{2}}{\ltos}
&\cH^{3,2}(2)&\tos \,\,   \cdots 
\stackrel{{\tilde{D}}_{n-3}}{\ltos}
&\cH^{n-2,2}(2)
  \\
& & \eqvert &&\congvert
&&\congvert && \congvert
&&\congvert&&\congvert
\\
0&\tos&\mathfrak{g}&\tos&
\cE^{1}(2)&\stackrel{D_{0}}{\ltos} &  
  \cE^{1,1}_0(2)&\stackrel{D_{1}}{\ltos}&
\cE^{2,2}_0(2)&\stackrel{D_{2}}{\ltos}
&\cE^{3,2}_0(2)&\tos\,\,   \cdots  
\stackrel{D_{n-3}}{\ltos}
&\cE^{n-2,2}_0(2) 

\end{array}
$$

\medskip
\noindent
where the differential operators are:
$$
\begin{aligned}
{\tilde{D}}_0&=\text{Sym}\,{\tilde{d}}_2\quad 
\text{(so }({\tilde{D}}_0\ft)_{IJ}=\partial_{(I}\ft_{J)})\\  
{\tilde{D}}_1&={\tilde{d}}_2{\tilde{d}}_1\\
{\tilde{D}}_r&={\tilde{d}}_1 
\quad r=2,3,\dots, n-3.  
\end{aligned}
$$ 
When $n=3$, the shortened complexes terminate with the spaces
$\cE^{1,1}_0(2)$, $\cH^{1,1}(2)$.  
\end{theorem}

Observe that the operators $\tilde{D}_r$ in the lifted complex are simpler
than their downstairs counterparts:  they do not involve the trace-free
part.  
For example, $\tilde{D}_0$ is the Killing operator, while $D_0$ is the
conformal Killing operator.  The conditions defining the
$\cH^{r,s}(w)$  
imply that the images of the $\tilde{D}_r$ are
already contained in trace-free tensors.   

As indicated above, the operator $D_{n-1}=\delta_2$ obstructs ambient lifts
of the next space $\cE^{n-1,1}_0$ in the deformation complex.  However, 
$\operatorname{im}D_{n-2}\subset \ker D_{n-1}$, and   
a section of $\cE^{n-1,1}_0$ which is in $\ker\delta_2$ does have 
an ambient lift to $\cH^{n-1,1}(0)$.  This lift is not unique.
Nonetheless, 
by appropriately modifying the lifted space, one can arrange a unique
ambient lift.  Thus it is possible to extend the above complexes one more
term to include an ambient realization of $\ker D_{n-1}$.  
For this term, the analogue of the restriction
operator inverse to the lift effectively involves a differentiation   
and the operator lifting $D_{n-2}$ has order one less than $D_{n-2}$.  When
$n>3$, the operator lifting $D_{n-2}:\cE^{n-2,2}_0(2)\to \ker D_{n-1}$ is 
$\tilde{d}_1:\cH^{n-2,2}(2)\to \cH^{n-1,2}(2)$.   When $n=3$, the
operator lifting  $D_1:\cE^{1,1}_0(2)\to \ker D_2$ is the second  
order operator 
$\tilde{d}_2\tilde{d}_1:\cH^{1,1}(2)\to \cH^{2,2}(2)$, 
the same operator which lifts $D_1$ in higher dimensions.    

Theorem~\ref{amlift} also gives ambient realizations for other  
gBGG complexes;  for example the deRham complex and the complex which   
resolves the standard representation $\R^{n+2}$ of $\mathfrak{g}$.    

Next we indicate how Theorem~\ref{jetisoodd} can be proved  
using the ambient lift of the deformation complex.    
In the previous section we constructed the map $c:\cM/\CDiff_0\to \cRt$  
which evaluates the curvature tensors of the ambient metric and outlined
why it is $P$-equivariant.  So what remains is to show that $c$ is 
bijective with polynomial inverse.  The first step is a linearization
argument as in the direct proof mentioned in \S\ref{odd}.  One truncates
all the jet spaces and 
the map $c$ at finite order to make everything finite-dimensional.
Geodesic normal coordinates and the ``conformal normal form'' mentioned
previously provide a slice for the $\CDiff_0$ action, from which it follows 
that $\cM^N/\CDiff_0$ is a smooth manifold, where $\cM^N$ indicates the
truncation of $\cM$ at order $N$.  Now either an algebraic induction
argument or the inverse function theorem can be used to reduce the
conclusion to proving that $dc:T\cM/T\cO\to T\cRt$ is a vector space
isomorphism, where $\cO$ is the $\CDiff_0$-orbit of the flat metric $h$,
$T\cM$ and $T\cO$ denote the tangent spaces at $h$, and 
$T\cRt$ is the tangent space to $\cRt$ at $0$.   

The second step is to relate the spaces $T\cM/T\cO$ and $T\cRt$ 
to the spaces appearing in the deformation complex.  
\begin{lemma}\label{iso1}
$T\cM/T\cO\cong \cJ^{1,1}_0(2)/D_0\cJ^1(2)$.  
\end{lemma}
\begin{proof}
The definitions give
\[
\begin{split}
T\cM &= \{s\in \cJ^{1,1}(2):s(0)=0\},\qquad \\
T\cO &= \{\mathcal{L}_Vh:V=O\left(\|x\|^2\right)\}\oplus  
\{\Omega^2h:\Omega = O\left(\|x\|^2\right)\}.
\end{split}
\]
Recall that $D_0$ is the conformal Killing operator, and corresponds to
$V\to\tf \mathcal{L}_Vh$ when its argument is viewed as a vector field.  
Now, as  
in the construction of geodesic normal coordinates, every 1-jet of 
an infinitesimal metric is in the range of the Killing operator on jets of
vector 
fields.  This shows that all 1-jets in both $T\cM$ and $\cJ^{1,1}_0(2)$ are 
contained in the respective denominator spaces.  For higher order jets,
the term $\{\Omega^2h:\Omega = O\left(\|x\|^2\right)\}$ cancels the   
trace components.
\end{proof}

\begin{proposition}\label{iso2}
$T\cRt\cong \ker\tilde{d}_1\subset \cJt{}^{2,2}_\cH(2)$.    
\end{proposition}
\begin{proof}
The tangent space $T\cRt$ is defined by the same relations (1)--(5) in 
Definition~\ref{Rt}, except that the $\tilde{Q}^{(s)}$ term in (3) is
replaced by 0.   Thus each $\Rt_{IJKL,M_1\cdots M_r}$ is symmetric in
$M_1\cdots M_r$.    
We may identify jets $\Rt$ at $e_0$ of sections of 
$\wt{}^2\otimes\wt{}^2$ with such lists of tensors by the requirement that  
\begin{equation}\label{jetident}
\pa^r_{M_1\cdots M_r}\Rt_{IJKL}(e_0)
=\Rt_{IJKL,M_1\cdots M_r},\qquad r\geq 0.
\end{equation}
Clearly conditions (1) and (4) are equivalent to the statement that 
$\Rt$ is the jet of a section of $\wt{}^{2,2}_0$.  Differentiation of the
relation $X^L\Rt_{IJKL}=0$ and evaluating at $e_0$ shows that condition 
(5) is equivalent to the statement that $X\into \Rt=0$ to infinite order.  
Condition (2) holds if and only if $\Rt\in \ker \tilde{d}_2$.  Since these 
are all the relations defining $\cRt$, it follows that 
$\ker \tilde{d}_2|_{\cJt{}^{2,2}_\cH(2)}\subset T\cRt$.   Note that  
$\ker \tilde{d}_1=\ker
\tilde{d}_2$ on sections of $\wt{}^{2,2}$ by the
symmetries of curvature tensors.

Now $\cJt{}^{2,2}_\cH(2)$ is defined by the conditions considered in the  
previous paragraph together with the 
additional requirements that $\Rt$ be homogeneous of degree 2 as a jet and
that $\Dt\Rt =0$ and $\dt_1\Rt=0$ to infinite order.   The homogeneity 
statement is equivalent to 
\begin{equation}\label{homo}
\Rt_{IJKL,M_1\cdots M_r0}=(-2-r)\Rt_{IJKL,M_1\cdots M_r}.
\end{equation}
The symmetry of $\Rt^{(r+1)}$ in the differentiation indices and relation
(2) can be used to express the left
hand side as a sum of two terms in which the '0' index is before the
comma.  
Then applying (5) and then (2) again establishes \eqref{homo}.  The
relations $\dt_1\Rt=0$ and $\Dt\Rt=0$ follow similarly using (2) to move 
contracted derivative indices before the comma and then applying (4).   
Thus $T\cRt=\ker \tilde{d}_1|_{\cJt{}^{2,2}_\cH(2)}$ under the
identification  \eqref{jetident}. 
\end{proof}

Composing with the isomorphisms of Lemma~\ref{iso1} and
Proposition~\ref{iso2}, the jet isomorphism theorem reduces to the
statement that 
$$
dc:\cJ^{1,1}_0(2)/D_0\cJ^1(2)\to \ker \tilde{d}_1\subset \cJt^{2,2}_\cH(2)     
$$
is an isomorphism.  Suppose first that $n\geq 5$.  According to
Theorem~\ref{amlift}, the lift of the deformation complex on jets contains  
\def\tos{\!\!\to\!\!}
$$
\begin{array}{ccccccc}
{\tos}&\cJt{}^{1}_\cH(2)&{\tos}&
  \cJt{}^{1,1}_\cH(2)&{\tos}&
 \cJt{}^{2,2}_\cH(2)&\stackrel{\tilde{D}_{2}}{\longrightarrow}   
  \\
&\congvert & &
 \congvert& &
\congvert& 
\\
{\tos}&\cJ^{1}(2)&\stackrel{D_{0}}{\longrightarrow} &  
  \cJ^{1,1}_0(2)&\stackrel{D_{1}}{\longrightarrow}&
\cJ^{2,2}_0(2)&\stackrel{D_{2}}{\longrightarrow}  
\end{array}
$$
Since the deformation complex is
exact on jets, $D_1$ induces an isomorphism 
$$
\cJ^{1,1}_0(2)/D_0\cJ^1(2)\cong \ker D_2 \cong \ker \tilde{D}_2 =\ker 
\tilde{d}_1.
$$
One can show that this map agrees with $dc$, and the result follows.

When $n=3$, $D_1$ maps into $\ker D_2\subset \cJ^{2,1}_0$.  But, as
discussed after the statement of Theorem~\ref{deflift}, when $n=3$ we have
a modified lift of  
$\ker D_2$ to $\cJt{}^{2,2}_\cH(2)$.    
The jet isomorphism theorem follows in exactly the same manner.

\section{Jet isomorphism, even dimensions}\label{even}

When $n$ is even, the construction of the ambient metric is obstructed at
order $n/2$.  So the map $c$ evaluating the covariant derivatives of
curvature of the ambient metric is not defined beyond this order.  This 
is a reflection of a difference in the structure of $\cM/\CDiff_0$ as a
$P$-space when $n$ is even.

The same phenomenon occurs when constructing the ambient lift for
$\cE^{r,s}_0(w)$ when $w+n/2-r-s\in \N$.  In this section, an  
extension of the theory to these cases will be outlined.  The main
ingredients are the following:
\begin{itemize}
\item
A weakening of the homogeneity condition on the ambient lift 
\item
The occurrence of logarithm terms in the solutions of the ambient equations   
\item
Existence of an ambiguity (nonuniqueness) in the solutions    
\item  
An invariant smooth part for the solutions with log terms   
\item
Jet isomorphism theorem for an enlarged space  
\end{itemize}
We will first illustrate the ideas by discussing the ambient lift with log
term and the generalization 
of Theorem~\ref{amlift} for scalars in the obstructed case $w+n/2\in \N$.
Then, by analogy with the discussion in \S\ref{odd}, 
we will formulate the jet isomorphism theorem  
for conformal structures in even dimensions and will discuss the
construction of inhomogeneous ambient metrics containing log terms and  
the extension of the map $c$ to infinite order in the even dimensional
case.  This is all joint work with K. Hirachi.  

Recall that Theorem~\ref{amlift} asserts that if $w+n/2\notin \N$, then
$\cE(w)\cong \cH(w)$, where $\cH(w)$ 
is the sheaf of harmonic jets along $\cN$ homogeneous of degree $w$.  But 
if $w+n/2=m\in\N$, then harmonic extension is obstructed at order $m$ by
the conformally invariant operator $\Delta^m$.  The following proposition 
shows 
that it is always possible to find a harmonic extension by including a log
term in the expansion.
\begin{proposition}\label{harmlog}
Suppose $w+n/2 =m \in \mathbb{N}$.  Let $\cU\subset \cQ$ be open and let 
$f\in \cE_\cU(w)$.  There exists an infinite order jet along
$\pi^{-1}(\cU)$ of a function $\ft$ on $\R^{n+2}$ of the form 
\begin{equation}\label{ftlog}
\ft = \st + \lt\, Q^m \log |Q|
\end{equation}
with $\st$, $\lt$ smooth, $\st$ homogeneous of degree $w$, 
$\lt$ homogeneous of degree $w-2m$, such that $\Dt \ft=0$ to infinite
order and 
$\ft|_{\pi^{-1}(\cU)} = f$.  These conditions uniquely determine $\lt$ to
infinite order along $\pi^{-1}(\cU)$ and determine $\st$ modulo 
$Q^m\cH_{\cU}(w-2m)$.       
\end{proposition}

Note that $\st$ and $\lt\, Q^m$ are each homogeneous of degree $w$.  Thus  
$\ft$ is almost homogeneous of degree $w$, but is not so because of the 
appearance of $\log |Q|$.  In this sense the homogeneity condition on $\ft$
has been weakened.  Of course, the appearance of this log term also means
that $\ft$ is no longer smooth.  

A main feature of
Proposition~\ref{harmlog} is that the solution $\ft$ is no longer unique.
One constructs $\ft$ inductively by order as in the proof of the scalar
case of Theorem~\ref{amlift} sketched in \S\ref{deform}.  The argument
there constructed $\ft \mod Q^m$.  The inclusion of the $Q^m\log |Q|$ term
enables 
the possibility of finding a harmonic extension at order $Q^m$.  The
coefficient of $Q^m\log |Q|$ is uniquely determined but not the coefficient
of $Q^m$, which can be prescribed arbitrarily on $\cN$.  The solution is
then uniquely determined to all higher orders.  The fact that the
uniqueness 
is at best modulo $Q^m\cH(w-2m)$ is immediate from \eqref{comm}:
$[\Dt,Q^m]=0$ on functions homogeneous of degree  
$w-2m=-n/2-m$.  Note that Theorem~\ref{amlift} implies that  
$\cH(-n/2-m)\cong \cE(-n/2-m)$ so that uniqueness modulo $Q^m\cH(-n/2-m)$
is  
the same as saying that the coefficient of $Q^m$ in the expansion of $\ft$
is undetermined.  This nonuniqueness is called the ambiguity in the
solution.     

It turns out that $\lt$ can be written entirely in terms of $\st$, and the
condition that $\ft$ be harmonic can be written entirely in terms of
$\st$.  Thus one can reformulate the extension as a map taking values in
a space of jets along $\cN$ of smooth homogeneous functions of degree
$w$, staying entirely in the smooth category.  To see this, straightforward
calculation using \eqref{comm} shows that 
\[
\begin{split}
\Dt\ft &= \Dt (\st + \lt\, Q^m \log |Q|)\\
&=(\Dt \st +4m\lt\, Q^{m-1})+(\Dt\lt)\,Q^m\log |Q|.  
\end{split}
\]
So $\Dt \ft=0$ to infinite order if and only if $\Dt \lt =0 \text{ and } 
\Dt \st =-4m\lt\, Q^{m-1}$ to infinite order.  Now iterating \eqref{comm}
shows that if $\lt$ is homogeneous of degree $-n/2-m$ and $\Dt\lt=0$, then 
\begin{equation}\label{cm}
c_m\Dt^{m-1}(Q^{m-1}\lt)=\lt,\qquad c_m^{-1}=(-4)^{m-1}(m-1)!^2.
\end{equation} 
Thus applying $c_m\Dt^{m-1}$ to the second equation gives 
$$
c_m \Dt^m \st = -4m\lt.
$$   
This gives $\lt$ in terms of $\st$.  Substituting back, one can write both  
equations in terms of $\st$:  
\begin{equation}\label{Hs}
\Dt\st = c_mQ^{m-1}\Dt^m\st \quad\text{ and }\quad \Dt^{m+1}\st=0.
\end{equation}
This motivates the following definition.
\begin{definition}
Suppose $w+n/2=m\in \N$.  Define $\cH_S(w)$ to be the sheaf on $\cQ$ of
infinite order jets along $\cN$ of smooth 
functions $\st$ on $\R^{n+2}$ which are homogeneous of degree $w$ and 
which satisfy \eqref{Hs} 
to infinite order along $\cN$, with $c_m$ as in \eqref{cm}.
\end{definition}

The conditions \eqref{Hs} are clearly $G$-invariant, so $\cH_S(w)$ is a
homogeneous sheaf on $\cQ=G/P$.   
Also observe by applying $c_m\Dt^{m-1}$ that if $c$ is any constant other
than 
$c_m$, then any solution $\st$ to the system obtained by replacing $c_m$ by
$c$ in \eqref{Hs} which is homogeneous 
of degree $w$ necessarily satisfies $\Dt^m\st=0$, and therefore $\Dt\st=0$.
The choice $c=c_m$ 
is the unique choice for which $\cH_S(w)\neq \cH(w)$.  

Now the substitute ambient lift theorem for scalars in the obstructed cases 
takes the form: 
\begin{theorem}\label{ambscalar}
Suppose $w+n/2 =m \in \mathbb{N}$.  There is a $G$-equivariant exact
sequence of sheaves:
\begin{equation}\label{exact}
0\to \cE(w-2m)\to \cH_S(w) \to \cE(w)\to 0.
\end{equation}
\end{theorem}

The map $\cH_S(w) \to \cE(w)$ is restriction to $\cN$.  
The map $\cE(w-2m)\to \cH_S(w)$ is harmonic extension followed by 
multiplication by $Q^m$; we saw in Proposition~\ref{harmlog} that 
jets in $Q^m\cH(w-2m)$ are already harmonic, so certainly are
contained in $\cH_S(w)$.  These maps are clearly $G$-equivariant.  
Since the sheaves $\cE(w)$ are soft, exactness of the sequence of sheaves 
is equivalent to  
exactness of the corresponding sequences of sections on any open set.     
And this is just the uniqueness statement of 
Proposition~\ref{harmlog} reformulated in terms of $\cH_S(w)$ as explained
above.   

By choosing a (necessarily non-$G$-equivariant) splitting of \eqref{exact},
one can parametrize $\cH_S(w)$ as $\cE(w)\times 
\cE(w-2m)$.  The space 
$\cE(w)$ corresponds to the initial density and $\cE(w-2m)$ to 
the ambiguity in the lift.  The space that has the ambient realization is
not the 
initial space $\cE(w)$ in which we were interested, but the enlargement
$\cE(w)\times \cE(w-2m)$ of this space by the ambiguity in the
solution.   
The space $\cH_S(w)$ realizing the ambient representation is an  
enlargement of the space of smooth homogeneous harmonic jets, and consists 
of the smooth homogeneous jets satisfying the system \eqref{Hs} rather than 
the equation $\Dt\ft=0$.  By taking jets at $e_0$ of the solutions of this
system, one obtains the substitute jet isomorphism theorem for scalars in
the obstructed cases analogous to the statement $\cJ(w)\cong \cJt_\cH(w)$
in the unobstructed cases.    

It is easy to check that when $n$ is even, all the spaces in the first half
of the deformation complex have $w+n/2-r-s\in \N$, 
so their ambient lifts are obstructed just as in the 
scalar case discussed above.  There is a
version of Theorem~\ref{ambscalar} for these spaces which is used in the
proof of the jet isomorphism theorem for conformal structures for $n$ even
as indicated below.

The jet isomorphism theorem for conformal structures in even dimensions
involves similar features as in the obstructed scalar case.  This time the
ambiguity is a symmetric 2-tensor which is trace-free with respect to the 
given metric.  
There is a map from jets of metrics together with jets of the ambiguity to
a space of ambient curvature tensors which induces a $P$-equivariant
bijection from the 
quotient by $\CDiff_0$.  We formulate the result more precisely.  

Define 
$$
\cM\times_\cM \cJ^{1,1}_0 \equiv \{(g,A)\in \cM\times \cJ^{1,1}:
g^{ij}A_{ij}=0 \text{  to infinite order}\}.   
$$
Here $\cJ^{1,1}$ denotes the space of jets of symmetric 2-tensors at $0\in
\R^n$ (ignoring the weight).  The  space $\cM\times_\cM \cJ^{1,1}_0$ may be 
regarded as a fiber bundle over $\cM$ by projecting onto the first factor.  
Set
$$
\cTt = \prod_{r=0}^\infty \textstyle{\bigwedge^{2,2}\R^{n+2}{}^* 
\otimes\bigotimes^r\R^{n+2}{}^*  \otimes \sigma_{-r-2}}.
$$
Here $\bigwedge^{2,2}\R^{n+2}{}^*$ denotes the finite-dimensional vector
space of 
covariant 4-tensors in $n+2$ dimensions with curvature tensor symmetries.    
Then $\cTt$ has a natural $P$-action. 
Recall that when $n$ was odd, the space $\cRt$ of lists of ambient 
curvature tensors was a $P$-invariant subset of $\cTt$.  The conformal jet
isomorphism theorem for $n$ even is then the following.

\begin{theorem}\label{jetisoeven}
Let $n\geq 4$ be even.  There is a $P$-equivariant polynomial injection    
$c:(\cM\times_\cM \cJ^{1,1}_0)/\CDiff_0\to \cTt$, whose image $\cRt$ is a   
submanifold of $\cTt$ whose tangent space $T\cRt$ at $0$ is the space of
jets $\Rt\in \cJt{}^{2,2}(2)$ which are solutions to the following  
equations to infinite order at $e_0$:  
\begin{enumerate}
\item
$\Rt_{IJ[KL,M]} = 0$  
\item
$X\into \Rt=0$
\item
$\trt\,\Rt = c_{n/2} Q^{n/2-1}\Dt^{n/2-1}\trt\,\Rt$   
\item
$\Dt^{n/2}\trt\,\Rt =0$.
\end{enumerate} 
Also, $c^{-1}:\cRt\to (\cM\times_\cM \cJ^{1,1}_0)/\CDiff_0$ is polynomial.    
Here $(\trt\,\Rt)_{IJ}=\htt^{KL}\Rt_{IKJL}$ corresponds to the Ricci
tensor, and $\cJt{}^{2,2}(2)$ is identified with a $P$-submodule of 
$\cTt$ via \eqref{jetident}.  The constant $c_{n/2}$ is given in
\eqref{cm}.  
\end{theorem}

The formulation of Theorem~\ref{jetisoeven} requires some explanation.   
First, the statement that $\cRt$ is a submanifold of $\cTt$ is to be  
interpreted in terms of finite-order truncations of these spaces; the full
spaces are projective limits of their truncations.  The truncations are
finite dimensional so these notions make sense in this context.  
Next, we have not yet defined the $\CDiff$-action on $\cM\times_\cM   
\cJ^{1,1}_0$ which gives rise to the quotient by $\CDiff_0$ and the
$P$-action on the quotient.  There is a natural action of $\CDiff$, but as
in Theorem~\ref{ambscalar}, it is not a 
product action.  This action has the property that the projection 
$\cM\times_\cM \cJ^{1,1}_0\to \cM$     
is $\CDiff$-equivariant, where the action on $\cM$ is that defined in
\S\ref{odd}.   
The $\CDiff$-action on $\cM\times_\cM \cJ^{1,1}_0$ will be defined below.     

When $n$ was odd, the nonlinear space $\cRt$ was identified explicitly; see  
Definition~\ref{Rt}.  For $n$ even, Theorem~\ref{jetisoeven} 
asserts instead that $c$ is a bijection onto a submanifold $\cRt$ of $\cTt$
and identifies explicitly the  
tangent space $T\cRt$.  This suffices for all the 
applications; the explicit form of the nonlinear terms is not needed.  
One can be somewhat more explicit about the equations defining $\cRt$.  
The equations (1)--(3) and (5) in Definition~\ref{Rt} hold also for $\cRt$
in even dimensions.  Equation (4) in Definition~\ref{Rt} is replaced by
nonlinear versions of (3) and (4) above.    

As shown by the proof of Proposition~\ref{iso2}, for $n$ odd $T\cRt$ is
defined by 
exactly the same relations (1)--(4) above, except that (3) and (4) 
are replaced by the single equation $\trt\,\Rt=0$.    
The relations (3) and (4) in Theorem~\ref{jetisoeven} are completely
analogous to the equations \eqref{Hs} for the obstructed 
scalar problem; the Ricci tensor $\trt\,\Rt$ plays the role of $\Dt \st$.  

In the rest of this section we will describe the extension of the ambient
metric construction to all orders in even dimensions and the construction
of the map $c$.  
Recall that in odd dimensions an ambient metric is a smooth metric defined
by the conditions (1)--(3) of Definition~\ref{gt}.  In even dimensions
there is a formal obstruction at order $n/2$ to the existence of such a
metric, analogous to the obstruction to finding a smooth harmonic extension
of a density in the scalar problem.  It is natural to try to continue the
expansion by including log terms.  In the scalar problem, $Q$ was a
$G$-invariant defining function for $\cN$ and terms involving $Q^m\log |Q|$ 
contained a built-in $G$-invariance.  But there is no canonical
analogue of $Q$ for the nonlinear problem, so it is not clear what the
argument of the logarithm should be to obtain an invariant construction.  
Another distinction is that the nonlinearity of the Ricci
curvature operator will force the inclusion of powers of the logarithm as 
well.  

These considerations motivate the following definition.  Let $r$
denote an arbitrary smooth defining function for $\cG\subset \cGt$
homogeneous of degree 2.
\begin{definition}
Let $\cA_{log}$ denote the space of formal asymptotic expansions of metrics
of signature $(p+1,q+1)$ on $\cGt$ of the form
\begin{equation}\label{metricexpansion}
\gt \sim \gt^{(0)}+ \sum_{N\geq 1}\gt^{(N)}r (r^{n/2-1}\log
|r|)^N,     
\end{equation}
where $\gt^{(N)}$, $N\geq 0$, are smooth symmetric 2-tensor fields on
$\cGt$ satisfying $\delta_s^* \gt^{(0)}=s^2\gt^{(0)}$ and 
$\delta_s^* \gt^{(N)}=s^{(2-n)N}\gt^{(N)}$ for $N\geq 1$, and such that   
$\iota^*\gt = {\bf g_0}$.   
\end{definition}
\noindent
It is easy to see that the space $\cA_{log}$ is independent of the choice
of $r$; upon changing $r$, one obtains an
expansion of the same form but with different coefficients.  Also, the
space 
$\cA_{log}$ is invariant under pullback by smooth homogeneous
diffeomorphisms $\Phi$ of $\cGt$ satisfying $\Phi|_\cG = I$.  

Recall that ambient metrics in odd dimensions automatically had an
additional geometric property.  We call metrics having this property
straight:  
\begin{definition}
A metric $\gt\in \cA_{log}$ is {\it straight} if for each $p\in\cGt$, 
the dilation orbit $s\to \delta_sp$ is a geodesic for $\gt$.  (Since $\gt$
is only defined as an asymptotic expansion, this means that the geodesic
equations hold to infinite order along $\cG$.)  
\end{definition}

The ambient metrics involving log terms in even dimensions are then defined
as follows.  We call these inhomogeneous ambient metrics because the
occurence of the log terms means that the metrics are no longer
homogeneous.   

\begin{definition}\label{ambmetric}
Let $n\geq 4$ be even.  An inhomogeneous ambient metric for $(M,[g])$ is a
straight metric $\gt\in \cA_{log}$ satisfying $\Ric(\gt)=0$ formally to
infinite order.        
\end{definition}

The straightness 
condition is crucial in the inhomogeneous case because of the following
proposition. 

\begin{proposition}\label{Tsmooth}
Let $\gt\in \cA_{log}$ be straight.  Then $\gt(T,T)$ is a smooth defining 
function for $\cG$ homogeneous of degree 2.
\end{proposition}
\noindent
Recall that in the flat case, the vector field $X$ on $\R^{n+2}$ plays the
role of $T$ and satisfies $\htt(X,X)=Q$.  Thus  
$\gt(T,T)$ is a generalization of $Q$.  For general $\gt\in \cA_{log}$,  
$\gt(T,T)$ has an asymptotic expansion involving $\log |r|$, but  
Proposition~\ref{Tsmooth} asserts that if $\gt$ is 
straight, then $\gt(T,T)$ is actually smooth.      
The proof of Proposition~\ref{Tsmooth} is a straightforward analysis of the 
geodesic equations for the dilation orbits.

Let $\gt\in \cA_{log}$ be straight. Then $\gt(T,T)$ is a  
canonically determined smooth defining function for $\cG$ 
homogeneous of degree 2.  We may therefore take $r=\gt(T,T)$ in
\eqref{metricexpansion}.  The term $\gt^{(0)}$ appearing in the resulting
expansion is then a smooth metric uniquely determined by $\gt$
independently of any choices.  We call this metric $\gt^{(0)}$ the 
{\it smooth part} of $\gt$.  Observe that $\gt^{(0)}$ is homogeneous of
degree 2.  One checks that $\gt^{(0)}$ is also straight.  
If $\Phi$ is a smooth homogeneous diffeomorphism satisfying $\Phi |_\cG =  
I$ and $\gt\in\cA_{log}$ is straight, then $(\Phi^*\gt)^{(0)}= 
\Phi^*(\gt^{(0)})$.  

We extend Definition~\ref{normalform} to the inhomogeneous case:  a
straight metric $\gt\in \cA_{log}$ is said to be in normal form relative to
a metric $g$ in the conformal class if its smooth part $\gt^{(0)}$ is in
normal form relative to $g$.  If $\gt\in\cA_{log}$ is straight, then 
there is a 
smooth homogeneous diffeomorphism $\Phi$ uniquely determined to infinite
order at $\rho=0$ such that $\Phi |_\cG = I$ and  
such that $\Phi^*\gt$ is in normal form relative to $g$.  

The main theorem concerning the existence and uniqueness of inhomogeneous
ambient metrics is the following.
\begin{theorem}\label{inhomo}
Let $n\geq 4$ be even.  
Up to pullback by a smooth homogeneous diffeomorphism which restricts to
the identity on $\cG$, the inhomogeneous ambient metrics for $(M,[g])$ 
are parametrized by the choice of an arbitrary trace-free symmetric
2-tensor field (the ambiguity tensor) on $M$.    
\end{theorem}

We describe more concretely the parametrization of inhomogeneous ambient
metrics in terms of the ambiguity tensor.  Choose a representative metric
$g$ in the conformal class; we normalize the diffeomorphism invariance by
requiring that $\gt$ be in normal form relative to $g$.  Let 
$\gt^{(0)}$ be the smooth part of $\gt$ and let $\cGt\cong \R_+\times
M\times \R$ be the decomposition induced by the choice of $g$.
Consider the component of
$\gt^{(0)}$ obtained by restricting to vectors tangent to $M$ in the
decomposition $\R_+\times M\times \R$; by homogeneity this may be written 
$\gt^{(0)}_{ij}=t^2 g_{ij}^{(0)}(x,\rho)$, where $g_{ij}^{(0)}(x,\rho)$  
is a smooth 1-parameter family of metrics on $M$ with $g_{ij}^{(0)}(x,0)$
equal to the chosen metric $g_{ij}(x)$.  The ambiguity tensor of $\gt$
relative to $g$ is:
$$
A_{ij}=\tf\left((\pa_\rho)^{n/2}g_{ij}^{(0)}|_{\rho =0}\right).  
$$
Theorem~\ref{inhomo} asserts that for each representative metric $g$ 
and each choice of trace-free symmetric 2-tensor $A$, there is 
a unique inhomogeneous ambient metric $\gt$ in normal form relative to $g$
with ambiguity tensor $A$.   

The choice of $g$ and $A$ uniquely determine $\gt$ in normal form, and
therefore also the smooth part $\gt^{(0)}$.  As in the scalar case, the
inhomogeneous ambient metric serves as an intermediate tool used to
determine the smooth 
homogeneous  metric $\gt^{(0)}$.  However, because of the nonlinearity we 
do not have a   
simple way of writing directly the system of equations defining
$\gt^{(0)}$.  

The map $c$ is
now defined exactly as in the odd-dimensional case, using the smooth part
$\gt^{(0)}$ in place of $\gt$.  Given $(g,A)\in \cM\times_\cM \cJ^{1,1}_0$, 
choose tensors also denoted $g$ and $A$ in a neighborhood of $0\in \R^n$
with the prescribed Taylor expansions and such that  
$g^{ij}A_{ij}=0$ in the whole neighborhood.  According to
Theorem~\ref{inhomo}, there is a unique 
inhomogeneous ambient metric $\gt$ in normal form relative to $g$ with
ambiguity tensor $A$.  Define the tensors $\Rt^{(r)}$ to be the iterated
covariant derivatives of the curvature tensor of $\gt^{(0)}$ evaluated at
$t=1$, $x=0$, $\rho=0$.     
This gives a map $\cM\times_\cM \cJ^{1,1}_0\to \cTt$.   
Because these $\Rt^{(r)}$ are the curvature tensors of a smooth,
homogeneous, straight metric, the relations  
(1)--(3) and (5) in Definition~\ref{Rt} hold for these tensors.  Since
$\gt^{(0)}$  is not Ricci-flat, equation (4) in Definition~\ref{Rt} does
not hold.  But a study of 
the linearized problem shows that relations (3) and (4) in
Theorem~\ref{jetisoeven} hold for the linearized tensors.  

The $\CDiff$ action on $\cM\times_\cM\cJ^{1,1}_0$ is defined as follows. 
First consider metrics $g$ and ambiguity tensors $A$ defined on a manifold
$M$; 
the $\CDiff$ action will be obtained by passing to jets at the origin in
$\R^n$.  
Let $\gt$ be the inhomogeneous ambient metric in normal form relative to
$g$ with ambiguity tensor $A$.  If $0<\Omega\in C^\infty(M)$, set $\gh =
\Omega^2g$.  Now there is a smooth homogeneous diffeomorphism $\Phi$
satisfying $\Phi|_\cG =I$, uniquely 
determined to infinite order at $\rho=0$, such that $\Phi^*\gt$ is in
normal form relative to $\gh$.  Since $\Phi^*\gt$ is also an inhomogeneous 
ambient metric, it uniquely determines an ambiguity tensor $\Ah$ with the
property that  
$\Phi^*\gt$ is the inhomogeneous ambient metric in normal form relative to
$\gh$ with ambiguity tensor $\Ah$.  The correspondence     
$(g,A,\Omega)\to \Ah$ gives a well-defined transformation law for the
ambiguity tensor under conformal change, described more explicitly in
\cite{GH1}.  The jet of $\Ah$ at a point depends only on the jets of 
$(g,A,\Omega)$ at that point.  Now for $(\varphi,\Omega)\in\CDiff$ and  
$(g,A)\in\cM\times_\cM\cJ^{1,1}_0$, the $\CDiff$ action is defined by
$$
(\varphi,\Omega).(g,A)
=\left((\varphi^{-1})^*\gh,(\varphi^{-1})^*\Ah\right). 
$$
One can identify the Jacobian along $\cG$ of the diffeomorphism $\Phi$
above  to derive the transformation laws for the tensors $\Rt^{(r)}$.  
From this it follows that the map $\cM\times_\cM \cJ^{1,1}_0\to \cTt$
passes to a map $c:(\cM\times_\cM \cJ^{1,1}_0)/\CDiff_0\to \cTt$ which is    
$P$-equivariant, as claimed in Theorem~\ref{jetisoeven}.  

The completion of the proof of Theorem~\ref{jetisoeven} requires showing
that 
$$
dc:T\left((\cM\times_\cM \cJ^{1,1}_0)/\CDiff_0\right)\to T\cRt
$$ 
is a vector space isomorphism.   This uses the same idea as for $n$ odd:    
lift the deformation complex.  However, the algebra is substantially more
complicated, as there is an ambiguity for the lift of each term in the
first half of the complex.    

It is possible to extend
the parabolic invariant theory of \cite{BEG} to characterize scalar
$P$-invariants of $T\cRt$ for $n$ even.   Theorem~\ref{jetisoeven} 
then enables one to transfer 
the results to characterize scalar invariants of conformal structures in
even dimensions similarly to the arguments of \cite{FG2} in odd dimensions.
These results are described in \cite{GH1}; details will be forthcoming.

\end{document}